\documentclass[12pt]{amsart}
\usepackage[only,ninrm,elvrm,twlrm,sixrm,egtrm,tenrm]{rawfonts}
\usepackage{subfigure}
\usepackage{graphicx}
\usepackage{rotating}
\usepackage{epsf}
\usepackage[T1]{fontenc}
\usepackage[utf8]{inputenc}
\usepackage{amssymb}
\usepackage{amsmath}    
\def\R{I\kern -0,37 em R}
\def\P{I\kern -0,37 em P}
\def\Z{I\kern -0,37 em Z}

\title{ON THE EQUIVALENCE PROBLEM FOR GEOMETRIC STRUCTURES, I}

\author{Antonio Kumpera}
\address[Antonio Kumpera]{Campinas State University\\Campinas, SP, Brazil}
\email{antoniokumpera@hotmail.com}

\date{May 2014}

\begin{document}

\begin{abstract}
We discuss the local and global problems for the equivalence of geometric structures of an arbitrary order and, in later sections, attention is given to what really matters, namely the equivalence with respect to transformations belonging to a given pseudo-group of transformations. We first give attention to general prolongation spaces and thereafter insert the structures in their most appropriate ambient namely, as specific solutions of partial differential equations where the equivalence problem is then discussed. In the second part, we discuss applications of all this abstract nonsense and take considerable advantage in exploring Élie Cartan's magical trump called \textit{transformations et prolongements mériédriques} that somehow seem absent in present day geometry.
\end{abstract}

\maketitle

\section{Introduction}
Structures have been studied for a long time and they might even be retraced back to Archimedes in his effort to confront and, if possible, to repel the Roman fleet during the siege of Syracuse. The \textit{structure} of his parabolic incendiary mirrors would burn up the roman vessels with their beams focused on the sails. Syracuse was, after two long years, finally taken over on a cloudy day by roman general Marcellus and Archimedes killed by a roman soldier. More than two thousand years later, Élie Cartan undoubtedly brought the most significant contributions due to his magical skill in dealing with \textit{continuous finite and infinite groups of transformations} (\cite{Cartan1904}). Our account, given here, is hopefully a non-tedious though, at best, just a partial repetition of what the reader can certainly find hidden somewhere in Cartan's \textit{O\!\!Euvres} (\cite{Cartan1952}). For this, however, he must not only read french but, above all, have \textit{Medusa}'s penetrating and petrifying stare. A less frightening alternative is offered by Mona Lisa's \textit{regard énigmatique} which entails her to stare directly into your eyes wherever you should be standing in the \textit{Louvre's} Main Hall. The reader should also be very careful since, as many claim, Cartan's \textit{O\!\!Euvres} might uncover a Pandora box.

\vspace{2 mm}

\noindent
The local equivalence problem for structures (or almost-structures) treated systematically and placed in a context of ample 
generality can be retraced back to Sophus Lie following the appearance, in the Mathematische Annalen (\cite{Lie1884}), of his renowned mémoire \textit{Über Differentialinvarianten} (an english translation being available in \cite{Ackerman1976}) and much later of his masterpiece \textit{Verwertung des Gruppenbegriffes für Differentialgleichungten, I} (\cite{Lie1895}). Needless to say, a structure for Sophus Lie acquired the more visible nature of a "geometric object" not resembling at all, at least at first sight, to the crushing formalisms of our present days and his theory of differential invariants had the precise aim of devising the "best"(shortest, most accurate and of lowest degree) integration methods for differential equations, based on the properties (structure) of the invariance groups. He showed, for example, that Jacobi's \textit{last multiplier} method was the best possible due to the fact that the pseudo-group of all volume preserving transformations is simple! In this restricted context, the structures stand out as solutions of differential equations (\cite{Kumpera1999}). It should be stated however that whereas Lie was incredibly successful in dealing with equations of "finite type" (the solutions depending only upon a finite number of parameters) since he had a deep knowledge and understanding of "finite continuous groups", in the general context he was unable to go any further beyond examining a few specific examples involving the simple infinite (transitive) groups. He was aware of all the four classes of (complex) simple groups but rather uneasy whether these were the only ones, thus preventing him to take any benefits stemming from a systematic use of Jordan-Hölder resolutions (\cite{Lie1895}). Consequently, it remained to Cartan to develop the infinite dimensional theory with a \textit{touch} only accessible to the most illuminated (see the references \cite{Cartan1905, Cartan1908, Cartan1909, Cartan1921, Cartan1930, Cartan1938/39}).  

\vspace{2 mm}
\noindent
In this paper we concentrate mainly on certain ideas relating the differential invariants of a Lie pseudo-group with the formal and local equivalence problem for structures. Apart Lie's basic ideas, we also find here Cartan and
Ehresmann's most inspiring sources (\cite{Cartan1910, Cartan1932, Ehresmann1958}). 
To avoid long repetitions, the technical aspects adopted here stem entirely from \cite{Kumpera1971} and \cite{Kumpera1972} though we try to simplify as far as possible the notations as well as avoid any excessive abstractions (in Donald Spencer's words, \textit{the abstract nonsense}). Furthermore, perhaps overdue attention has been given to the differential structure of the equations defining $k-th$ order equivalences of structures. Whereas in the theory of finite dimensional Lie groups the three fundamental theorems of Lie shine as neat and as beautiful as Botticelli's \textit{La Nascita di Venere} and Ingre's \textit{Le Printemps}, in the case of groupoids and Lie pseudo-groups they look more like Picasso's \textit{Guernica} and are a painful headache. In our case however, we only have to cope with the second theorem that fortunately is, in the present situation, very condescending and indulgent (sect.4,5,6, \cite{Cartan1937}, \cite{Cartan1938}).

\vspace{2 mm}

\noindent
Last but not least, here are a few words for "peer reviewers". Curiously enough, Sophus Lie and Élie Cartan did always row off the "main stream" for the simple reason that, at their time, essentially nobody was able to understand their writings. It took Cartan to understand what Lie did mean and Charles Ehresmann to understand Cartan. Significantly enough, Lie was Cartan's thesis adviser and Cartan was Ehresmann's adviser. As for Sophus Lie, he in fact never needed any adviser at all since he began writing in Norwegian so nobody would understand him anyway. Unfortunately, the author is unable to pinpoint anybody who did or who does (or who ever will) \textit{really} understand the full extent of Ehresmann's thoughts in all their galactic magnitude. Most probably we shall have to await for the next millennium\footnote{Il est à remarquer cependant que Ehresmann, bien qu'il construisait de très beaux arcs-en-ciel, ne s'est jamais soucié d'aller chercher le trésor se trouvant au bout. Par contre, Lie ainsi que Cartan allaient chercher désormais ce trésor sans se soucier à peindre au préalable de beaux arcs-en-ciel.}.

\section{Ehresmann's prolongation spaces}
Let \textit{P} be a differentiable manifold where, for convenience, we assume all the data of class $C^\infty$ though it would suffice to assume differentiability just up to a certain order. A \textit{finite} prolongation space of \textit{P} is a quadruple $(E,\pi,P,p)$ where \textit{E} is a differentiable manifold called the total space of the prolongation, $\pi:E\longrightarrow P$ a fibration (surmersion) and \textit{p} a prolongation operator that associates to each local diffeomorphism $\varphi$ of \textit{P} a local diffeomorphism $p\varphi$ of \textit{E} whose source and target are $\pi$-saturated open sub-sets inverse images of the source and target, respectively, of $\varphi$ and that furthermore obey the following requirements:

\vspace{4 mm}

\hspace{6 mm}\textit{i}) $p\varphi$ commutes with $\varphi$ and the projection $\pi$,

\vspace{2 mm}

\hspace{5 mm}\textit{ii}) \textit{p} is local and preserves pastings (\textit{recollements}),

\vspace{2 mm}

\hspace{4 mm}\textit{iii}) \textit{p} is a groupoid functor with respect to local diffeomorphisms ($\varphi$ being composable with $\psi$ whenever $\alpha(\psi)\cap\beta(\varphi)$ is non void, the unities being the identities on the open sets),

\vspace{2 mm}

\hspace{5 mm}\textit{iv}) Every differentiable one-parameter family $(\varphi_t)$ of local diffeomorphisms of \textit{P} prolongs, by \textit{p}, onto a one-parameter family of local diffeomorphisms of \textit{E} for which the vector field $d/dt(p\varphi_t)_{t=0}$ depends only upon $d/dt(\varphi_t)_{t=0}$ and projects onto it by $T\pi$.

\vspace{4 mm}

\noindent
We shall say that $p\varphi$ is the prolongation of $\varphi$ and, in order to simplify notations, the prolongation space will just be denoted by \textit{E}.

\noindent
Much in the same way, an \textit{infinitesimal} prolongation space of \textit{P} is a quadruple $(E,\pi,P,p)$ where the prolongation operator \textit{p} associates to each local vector field (infinitesimal transformation) $\xi$ given on \textit{P}, a \textit{prolonged} vector field $p\xi$ defined on the inverse image of the source $\alpha(\xi)$, this operation satisfying the corresponding (infinitesimal) properties:

\vspace{4 mm}

\hspace{6 mm}\textit{i}) $p\xi$ is $\pi$-projectable onto $\xi$,

\vspace{2 mm}

\hspace{5 mm}\textit{ii}) \textit{p} is local and preserves pastings,

\vspace{2 mm}

\hspace{4 mm}\textit{iii}) \textit{p} is a pre-sheaf morphism of Lie algebras.

\vspace{4 mm}

\noindent
Any finite prolongation space determines uniquely an infinitesimal prolongation space by which it is generated but what really matters is the converse that is not always true as we shall see in the sequel. Most fibre bundles considered in geometry (\textit{e.g.}, tensor bundles, Cartesian frames and co-frames, Stiefel truncated frames and co-frames, Grassmannian contact elements and their corresponding higher order analogues) are of course finite or infinitesimal prolongation spaces or both though our main interest is directed towards jet spaces. Let us also observe that prolongation spaces "compose" since the prolongation algorithm itself can be composed.

\vspace{4 mm}

\noindent
Let $\pi_0:P \longrightarrow M$ be a fibration (surmersion), denote by $J_k P$ the $k-th$ order jet bundle of local sections of $\pi_0$ and $\alpha$, $\beta$, $\rho_{hk}$ the well known projections. Following Ehresmann, we also denote by $\Pi_k P$ the groupoid of all invertible $k-$jets of the manifold \textit{P} ($k-$jets of local diffeomorphisms), by $J_k TP$ the vector bundle of all $k-$jets of local sections of $~TP\longrightarrow P~$ \textit{i.e.}, the jets of local vector fields on \textit{P} and finally by $\tilde{J}_k TP$ the fibration of all $k-$jets of local sections of the composite fibration $TP\longrightarrow P\longrightarrow M$. In the sequel, this \textit{tilde} notation will always be used for jets of sections of composite fibrations.

\newtheorem{prol}[LemmaCounter]{Lemma}
\begin{prol}
Given any fibration $\pi_0:P\longrightarrow M$, the jet space $J_k P$ has a natural infinitesimal prolongation space structure $(J_k P,\beta,P,p_k)$ where $\beta$ is the target map and where $p_k$ is the $k-th$ order standard  prolongation operator for  vector fields (by the target).
\end{prol}

\vspace{4 mm}

\noindent
Though the prolongation morphism $p_k$ goes back to Sophus Lie (in coordinates), we would like to add a few words so as to avoid any misunderstanding. It is \textit{not} possible to prolong, to $J_k P$, any local diffeomorphism $\varphi$ of \textit{P} since such a map can upset transversality (generic position) of a local section with respect to the $\pi_0$-fibres. However, when this condition is fulfilled, we can transform (at least locally) any $\pi_0$-section, whose image is contained in the domain of $\varphi$, into a new $\pi_0$-section and thereafter take its \textit{k}-jet. In particular, we shall then be able to prolong any $\pi_0$-projectable local diffeomorphism $\varphi$ of \textit{P} and the prolongation functor thus obtained, on projectable maps only, will of course fulfill the above stated properties of a finite prolongation space. Since we can always define a local vector field by its (local) one-parameter group $(\varphi_t)_t$ and since $\varphi_0=Id$, there is no restriction whatsoever in the prolongation procedure, to $k-th$ order, of \textit{any} local vector field defined on \textit{P} whereupon $J_k P$ becomes an authentic \textit{infinitesimal} prolongation space of \textit{P}.

\vspace{2 mm}

\noindent
Inasmuch, we can say that $TP,~J_k TP,~TM~and~J_k TM$ are prolongation spaces (finite and infinitesimal) of \textit{P} and \textit{M} respectively and, moreover, that $\tilde{J}_k TP$ is an infinitesimal prolongation space not only of \textit{TP} but also of \textit{P} for we can first prolong the local vector field $\xi$, defined on \textit{P}, to \textit{TP} and thereafter proceed with the above described "jet space prolongation". It should also be noted that $\tilde{J}_k TP$ is a (locally trivial) vector bundle with base space \textit{P} since the \textit{k-th} order tangency renders possible the vector space operations on the fibres. Finally, we would like to "stress" the condition (\textit{iii}) above by writing explicitly the equality $p(f\xi)=fp\xi$ where \textit{f} is any function.

\vspace{2 mm}

\noindent
As is well known, $(h+k)-$jets can become $h-$jets of $k-$jets and, inasmuch, $(h+k)-$jets can operate on "split" jets this motivating the following definitions:

\newtheorem{finite}[DefinitionCounter]{Definition}
\begin{finite}
a) The finite prolongation space E is said to be of order $\ell$ when, for any $k\geq 0~$, the k-jet of $p\varphi$ at the point $z\in E$ only depends upon the $(\ell+k)-$jet of $\varphi$ at the point $y=\pi (z)\in P$.

\vspace{2 mm}

\hspace{24 mm}b) The infinitesimal prolongation space E is said to be of order $\ell$ when, for any $k\geq 0$, the $k-$jet of $p\xi$, at the point $z\in E$, only depends upon the $(\ell+k)-$jet of $\xi$ at the point $y=\pi(z)\in P$.
\end{finite}

\vspace{2 mm}

\noindent
Recalling that the jet bundle $J_k TM$ identifies with the Lie algebroid of the Lie groupoid $\Pi_k M$, we can define for the above prolongation spaces of finite order and for any fixed positive integer \textit{k}:

\vspace{2 mm}

a) A left action

\begin{equation}
\Lambda_k:\Pi_{\ell+k} P~\times_P~\Pi_kE \longrightarrow \Pi_kE
\end{equation}

\vspace{2 mm}

\noindent
of the Lie groupoid $\Pi_{\ell+k}P$ on the groupoid $\Pi_{k}E$ (one can actually replace the last groupoid by the space of all $k-$jets $J_k(E,E)$)  by setting

\begin{equation}
(j_{\ell+k}\varphi(\beta X),X) \longmapsto j_k(p\varphi)(\beta X)\cdotp X~,
\end{equation}

\vspace{2 mm}

\noindent
the fibre product being taken with respect to $\alpha$ , for the first factor, and with respect to $\pi\circ\beta$ , for the second factor,

\vspace{2 mm}

b) its infinitesimal generator namely, the morphism

\begin{equation}
\lambda_k:J_{\ell+k}TP~\times_P~E \longrightarrow J_k TE
\end{equation}

\vspace{2 mm}

\noindent
defined by

\begin{equation}
(j_{\ell+k}\xi (y),X)\longmapsto j_k(p\xi)(X)~,~y=\pi X~,
\end{equation}

\vspace{2 mm}

c) and, finally, its extension to an infinitesimal action on the jets of the tangent bundles

\begin{equation}
T\Lambda_k:J_{\ell+k+1}TP~\times_P~J_{k+1}TE \longrightarrow J_k TE
\end{equation}

\vspace{2 mm}

\noindent
defined by the left Lie bracket action

\begin{equation}
(j_{\ell+k+1}\xi(y),X)\longmapsto [j_{k+1}(p\xi)(\alpha X),X]~,~y=\pi\alpha X~,
\end{equation}

\vspace{2 mm}

\noindent
where we are forced to augment the order by 1 since brackets absorb one order of differentiation. The reason for putting in evidence the Lie bracket becomes apparent if we operate, as is usually done, by taking local one parameter groups and thereafter differentiating.

\vspace{2 mm}

\noindent
We now claim that the action $\Lambda_k$ is differentiable. In fact, using standard methods involving the Lie algebroid $J_kTP$ of the groupoid $\Pi_k P$ (or, if one prefers, the sheaf $\underline{J_kTP}$), we can define a local exponential map that, with the help of the property (\textit{iv}) will provide the required differentiability. It then also follows that the infinitesimal generator $\lambda_k$ as well as the infinitesimal action $T\Lambda_k$ are differentiable though this property can be proved directly by observing that the vector bundles involved are locally trivial and generated by local holonomic  sections. It should also be observed that $T\Lambda_k$ is bilinear over \textbf{R} when the source and target spaces are fibered over \textit{E}.

\vspace{2 mm}

\noindent
We next observe that given a finite number of "composable" prolongation spaces, each of finite order, the composed prolongation space is also of finite order equal to the sum of the individually prescribed orders.

\vspace{2 mm}

\noindent
Other prolongation spaces that will be of our interest are those described in the following

\newtheorem{higher}[LemmaCounter]{Lemma}
\begin{higher}
To each finite or infinitesimal prolongation space $(E,\pi,P,p)$ and to each positive integer h corresponds, in a canonical way, a prolongation space $(J_h E,\alpha,P,p_h)$ verifying the following properties:

\vspace{2 mm}

a) When E is of finite order $\ell$ then $J_h E$ is of order $\ell+h$.

\vspace{2 mm}

b) The target projection $\beta:J_h E\longmapsto E$ is a surjective morphism of prolongation spaces i.e., respects the fibrations over P and commutes with the respective prolongation operations.

\vspace{2 mm}

c) More generally, the projection $\rho_{k,h}:J_h E\longrightarrow J_k E$ is a prolongation spaces morphism.
\end{higher}

\vspace{2 mm}

\noindent
Here, $J_h E$ is the set of all \textit{h}-jets of local sections of $\pi$ and is called the standard \textit{h-th order}  prolongation of \textit{E}, the prolongation operation being the composite of the operation provided by \textit{E} followed by the standard jet prolongation. Needless to say, everything that was stated concerning the prolongation space $J_k P$ in the Lemma 1 can be paraphrased \textit{ipsis litteris} for the above data. Inasmuch, we can also repeat everything that was said previously for the prolongation space $\tilde{J}_k TP$ of \textit{infinitesimal variations} relative to the composite 

\begin{equation*}
TP\longrightarrow P\longrightarrow M~,
\end{equation*}

\vspace{2 mm}

\noindent
the second arrow being equal to $\pi_0$, as well as for the $k-th$ order variations space $\tilde{J}_k E$ composed of all \textit{k}-jets of local sections of the composite fibration $E\longrightarrow P\longrightarrow M$, where $E\longrightarrow P$ is a prolongation space and $P\longrightarrow M$ simply a fibration giving rise to the finite or infinitesimal variations (\textit{cf.}, \cite{Kumpera1975} for the definitions). As for the bracket operation considered previously when defining an infinitesimal action, it is well defined in the present context due to the contact order conditions imposed. Finally, all the above considerations also extend to pre-sheaves $\Gamma(~)$ of local sections and enable us to operate with locally defined objects. In the sequel we shall also need the following extension of (3) namely,

\begin{equation}
\overline{\lambda}_k:J_{\ell+k}TP~\times_P~J_k E \longrightarrow TJ_k E~,
\end{equation}

\vspace{2 mm}

\noindent
defined by $(j_{\ell+k}\xi,j_k\sigma)= p_k\circ p(\xi)(j_k\sigma)$ and where $TP\longrightarrow P$ is the tangent prolongation space. In much the same way and using the prolongation operator, we can define the morphism

\begin{equation}
\tilde{\lambda}_k:J_{\ell+k}TP~\times_P~J_k E \longrightarrow J_k TE~,
\end{equation}

\vspace{2 mm}

\noindent
as well as the \textit{semi-holonomic} extension

\begin{equation}
\overline{\lambda}_{k+h}:J_{\ell+k+h}TP~\times_P~J_k E \longrightarrow J_h(TJ_k E)~.
\end{equation}

\vspace{2 mm}

\noindent
We thus see that the choices are many and, in fact, we could go on much further with the Ehresmannian game of the \textit{jeu de la théorie des jets} by considering for instance semi-holonomic, sesqui-holonomic and (definitely) non-holonomic jets but fortunately these will be of no purpose to us so we might as well forget about them right away. Furthermore, and this will be very useful, we can play the Ehresmannian game with \textit{differential forms} and \textit{co-tangent bundles} that, in this case, will act \textit{co-variantly}. As for the prolongation operation, there is of course essentially just one such operation that can however be vested under two or three garbs. Instead of prolonging by the target, as is done in the present paper, we can also prolong by the source or even, combining the two, we can prolong via the \textit{anchor}, the same one that holds Kirill anchored and not on the sail. In later sections we shall also introduce \textit{merihedric} prolongations (\textit{prolongements mériédriques de Élie Cartan}) as well as such transformations on jet spaces and higher order Grassmannians. Contrary to what is standard, there are uncountably many possibilities for the merihedric functor each one having its own merits and outstanding performance. In fact, Medusa as well as Mona Lisa claim that the merihedric setup was the magical trump and joker hidden in Cartan's sleave. Fortunately we shall need not talk about merihedric jet spaces, the standard ones being still of good use. One last remark: When differentiability is replaced by analyticity in the initial requirements for the prolongation spaces then, of course, all the other data also become analytic. We hope as well that the reader already noticed our small notational changes. Instead of the standard $j^k_x\sigma~$, we prefer to write $j_k\sigma(x)$ and, instead of $J^k$, we write  $J_k$, such notations rendering more pleasant and "co-variant" their composites.

\section{Symbols}
In this section we shall often refer to former publications (\cite{Kumpera1971}, \cite{Kumpera1972}) for the notations and the results which, however, are completely standard and well known. In this sense, we denote by \textit{D} the Spencer operator, by $S^k$ the bundle of symmetric tensors and by $\delta$ the Spencer operator restricted to the principal parts.
\noindent
Our first task is to examine what happens with the \textit{kernels} and consequently write:

\begin{equation*}
0\longrightarrow S^{\ell+k}T^*P\otimes TP\times_P E\longrightarrow J_{\ell+k}TP\times_P E\longrightarrow J_{\ell+k-1}TP\times_P E\longrightarrow 0 \end{equation*}
\begin{equation*}
\hspace{20 mm}\downarrow\ell_k\hspace{32 mm}\downarrow\lambda_k\hspace{27 mm}\downarrow\lambda_{k-1}\hspace{4 mm}
\end{equation*}
\begin{equation*}
0\longrightarrow\hspace{4 mm}S^k T^*E\otimes TE\hspace{8 mm}\longrightarrow\hspace{7 mm} J_k TE\hspace{7 mm}\longrightarrow\hspace{8 mm} J_{k-1}TE\hspace{4 mm}\longrightarrow 0
\end{equation*}

\vspace{4 mm}

\noindent
as well as

\begin{equation*}
0\!\rightarrow\!S^{\ell+k}T^*P\otimes TP\times_P J_k E\!\rightarrow\!J_{\ell+k}TP\times_P J_k E\!\rightarrow\!J_{\ell+k-1}TP\times_P J_{k-1} E\!\rightarrow\! 0       
\end{equation*}
\begin{equation*}
\hspace{12 mm}\downarrow\overline{\ell}_k\hspace{31 mm}\downarrow\overline{\lambda}_k\hspace{31 mm}\downarrow\overline{\lambda}_{k-1}
\end{equation*}
\begin{equation*}
0\rightarrow S^k T^*E\otimes VE\times_E J_kE\hspace{2 mm}\longrightarrow\hspace{2 mm} TJ_kE\hspace{3 mm}\longrightarrow\hspace{2 mm}TJ_{k-1} E\times_{J_{k-1}E} J_kE\!\rightarrow\! 0
\end{equation*}

\vspace{4 mm}

\noindent
and observe that the above $\overline{\lambda}_k$ is equal to $\lambda_0$ relative to the prolongation space $J_k E$ of \textit{P}, the later being of finite order $\ell+k$. More generally and with the obvious notations, $\overline{\lambda}_{k+h}(E)=\lambda_h(J_k E)$. Next, we claim that each family $(\ell_k)$ and $(\overline{\ell}_k)$ is a \textit{natural transformation} of the corresponding $\delta-cohomology~complexes$ this becoming apparent by examining the commutativity of the two diagrams below and observing that the vertical map $\overline{\ell}_k$ in the second diagram is in fact equal to

\begin{equation*}
[S^{\ell+k}T^*P\otimes TP\times_P J_1 E]\times_{J_1 E} J_k E\xrightarrow{\overline{\ell}_k\times Id}[S^k T^*E\otimes VE\times_E J_1 E]\times_{J_1 E} J_k E
\end{equation*}

\vspace{2 mm}

\noindent
the term $\overline{\ell}_k$ depending only upon its projection in $J_1 E$. In order to show the naturality we are forced, of course, to confront the above maps with the Spencer operator hence the necessity in extending the actions to the sheaf level. 

\begin{sidewaysfigure}[htp!]

\vspace{125 mm}

\begin{equation*}
0 \longrightarrow S^{\ell+k}T^*P~\otimes~TP~\times_P~E\xrightarrow{\delta\times Id}T^*P~\otimes~S^{\ell+k-1}T^*P~\otimes~TP~\times_P~E\xrightarrow{\delta\times Id}\wedge^2T^*P~\otimes~S^{\ell+k-2}T^*P~\otimes~TP~\times_P~E\hspace{1 mm}\xrightarrow{\delta\times Id}\hspace{1 mm}\bf{\cdots}
\end{equation*}

\begin{equation}
\hspace{9 mm}\downarrow\ell_k \hspace{51 mm}\downarrow\pi^*~\otimes~\ell_{k-1}\hspace{50 mm}\downarrow\pi^*~\otimes~\ell_{k-2}
\end{equation}

\begin{equation*}                
0\longrightarrow \hspace{4 mm}S^k T^*E~\otimes~TE\hspace{10 mm}\xrightarrow{\delta}\hspace{8 mm}T^* E~\otimes~S^{k-1} T^*E~\otimes~TE\hspace{6 mm}\xrightarrow{\delta}\hspace{8 mm}\wedge^2T^*E~\otimes~S^{k-2} T^*E~\otimes~TE\hspace{7 mm}\xrightarrow{\delta}\hspace{4 mm}\bf{\cdots}
\end{equation*}

\vspace{30 mm}

\begin{equation*}
0\rightarrow S^{\ell+k}T^*P\otimes TP\times_P J_1 E\xrightarrow{\delta\times Id}\hspace{1 mm}T^*P\otimes S^{\ell+k-1}T^*P\otimes TP\times_P J_1 E\xrightarrow{\delta\times Id}\hspace{1 mm}\wedge^2T^*P\otimes S^{\ell+k-2}T^*P\otimes TP\times_P J_1 E\xrightarrow{\delta\times Id}\bf{\cdots}
\end{equation*}

\begin{equation}
\hspace{1 mm}\downarrow\overline{\ell}_k\hspace{43 mm}\downarrow Id_{T^*P}~\otimes~\overline{\ell}_{k-1}\hspace{34 mm}\downarrow Id_{\wedge^2T^*P}~\otimes~\overline{\ell}_{k-2}
\end{equation}

\begin{equation*}
0\rightarrow\hspace{1 mm}S^kT^*E\otimes VE\times_E J_1 E\xrightarrow{\delta\times Id}T^*P\otimes S^{k-1}T^*E\otimes VE\times_E J_1 E\xrightarrow{\delta\times Id}\hspace{1 mm}\wedge^2T^*P\otimes S^{k-2}T^*E\otimes VE\times_E J_1 E\xrightarrow{\delta\times Id}\bf{\cdots}
\end{equation*}

\end{sidewaysfigure}

\noindent
We first consider the diagram (10), extend the infinitesimal generator (3) to the pre-sheaf of local sections and derive thereafter the sheaf morphism

\begin{equation*}
\underline{\lambda}_k:\underline{J_{\ell+k}TP}~\times_P~E~\longrightarrow~\underline{J_kTE},
\end{equation*}

\vspace{2 mm}

\noindent
where $\underline{\lambda}_k(\underline{\sigma}_y,q)=(\underline{\Gamma(\lambda_k)\sigma})_q~,~y=\pi(q)~,~\Gamma(~)$ denoting the set of all local sections. It now suffices to show that $\underline{\lambda}_k$ commutes with \textit{D} or, in other words, that the diagram (12) below commutes:

\vspace{8 mm}

\begin{equation*}
\underline{J_{\ell+k}TP}~\times_P~E~\xrightarrow{D~\times~Id}~(\underline{T^*P}~\times_P~E)~\otimes~
(\underline{J_{\ell+k-1}TP}~\times_P~E)
\end{equation*}
\begin{equation}
\downarrow\underline{\lambda}_k\hspace{30 mm}\downarrow~\underline{\pi}^*\otimes\underline{\lambda}_{k-1}
\end{equation}
\begin{equation*}
\underline{J_kTE}\hspace{6 mm}\xrightarrow{D}\hspace{10 mm}\underline{T^*E}~\otimes~\underline{J_{k-1}TE}\hspace{6 mm}
\end{equation*}

\vspace{4 mm}

\noindent
The commutativity is obvious on holonomic sections since \textit{D} vanishes on them and thereafter it suffices to argue as in \cite{Kumpera1972}, pg.75, using the formula (2) shown on that page. The restriction to the symbols then proves our claim for the first square of (10). To complete the proof, it suffices to adapt the diagram (12) to the other squares of (10) by iterating successively the Spencer operator. However, the previously mentioned formula (2) implies immediately the desired result as soon as the first square commutes. The commutativity of (11) is proved essentially in the same manner by observing that $\overline{\lambda}_k$ is the pullback of $\lambda_k$ via the composite $p_k\circ\sharp~$, where $~\sharp:J_kTE~\times_E~J_kE\longrightarrow \tilde{J}_kTE~$ is the fibre bundle morphism considered in \cite{Kumpera1971}, Propos.11.1 (\textit{M} standing in the  place of \textit{P}). The proof can then be achieved with the help of the \textit{reduced form of the holonomic prolongation operator} acting on infinitesimal variations and the Theorem 12.1.

\vspace{4 mm}

\noindent
Since the $(\ell_k)$ and $(\overline{\ell}_k)$ are natural transformations of the $\delta$ complexes and since the $(\lambda_k)$ and $(\overline{\lambda}_k)$ commute with \textit{D} it also follows that these natural transformations are compatible (commute) with the algebraic prolongations operating on the principal parts (i.e., on the symbols).

\vspace{4 mm}

\noindent
For the sake of not omitting any useful information, let us finally observe that the total spaces of the finite prolongation spaces are always locally trivial bundles  over their base spaces, a fact that is not true anymore for infinitesimal prolongation spaces (\textit{e.g.}, withdraw a point from the total space). However, when the fibres over the base space are compact and connected then these total spaces also become locally trivial, this last claim being just a special case of very general results due to Ehresmann.

\section{Lie and Cartan's notion of structure}
Let $(E,\pi,P,p)$ be a finite or infinitesimal prolongation space.

\newtheorem{struct}[DefinitionCounter]{Definition}
\begin{struct}
A finite or infinitesimal (almost-)structure of species E on the manifold \text{P} is the data provided by a global section (continuous, differentiable, analytic) of the prolongation space \textit{E}.
\end{struct}

\vspace{2 mm}

\noindent
We can define as well a structure of species \textit{E} above an open set \textit{U} of \textit{P} by simply taking a local section\footnote{Il serait malheureux de l'appeler une \textit{structure locale} car cette terminologie, due à Ehresmann, a un sens très précis autre que ci-dessus. Toute structure d'espèce \textit{E} est une espèce de structure locale.} and thereafter derive the notions of \textit{germ of structure, $k-$jet of structure, $k-th$ order contact element of a  structure} inasmuch as an \textit{infinite jet of a structure} (formal structure at a point). We should also mention that the above definitions are not any of Spencer's \textit{abstract nonsense} since many of the well known structures (\textit{e.g.}, Riemannian, conformal, almost-complex, almost-symplectic, etc.) are defined by global sections of locally trivial tensor bundles with fibres homogeneous spaces of linear groups. Inasmuch, almost-product, contact, Stiefel and many other structures are defined by sections of well known bundles in homogeneous spaces. Moreover, as we shall see further, the \textit{integrability} of an almost-structures can be detected by a Pfaffian system (also a structure!) defined on \textit{E} or on one of its prolongation spaces.

\vspace{2 mm}

\noindent
Given a structure \textit{S} of species \textit{E} (a finite prolongation space), associated to it are the pseudo-group $\Gamma(S)$ of all its local automorphisms as well as the infinitesimal pseudo-algebra $\mathcal{L}(S)$ of all its local infinitesimal automorphisms (vector fields). By definition, $\varphi\in\Gamma(S)$ when $p\varphi$ leaves invariant the sub-manifold \textit{im S} or, in other terms, when $\varphi(S)\stackrel{def}{=}p\varphi\circ S\circ\varphi^{-1}=S$ in the appropriate domains. Much in the same way, $\xi\in\mathcal{L}(S)$ when $\xi$ generates a local one-parameter group $(\varphi_t)$ with all its elements in $\Gamma(S)$ which amounts to say that $p\xi$ is tangent to \textit{im S} or equivalently that the Lie derivative $\theta(p\xi)S=\frac{d}{dt}~\varphi_t^{-1}(S)|_{t=0}$ vanishes (for each fixed $Y\in im~S$, take the tangent vector to the curve issued from \textit{Y}). Quite often, it is convenient to consider (differentiable) one-parameter families of local transformations that are not necessarily local one-parameter groups. With that in mind, let us next remark that such a local one-parameter family $(\varphi_t)$ defined in \textit{P} is entirely composed of elements belonging to $\Gamma(S)$ if and only if the following two properties are fulfilled:

\vspace{2 mm}
       
a) $\varphi_{t_0}\in\Gamma(S)$ for some value $t_0$ (\textit{t} varies in an open interval) and

\vspace{2 mm}

b) for each \textit{u}, $T\varphi_u^{-1}\circ\frac{d}{dt}~(\varphi_t)|_{t=u}\circ\varphi_u\in\mathcal{L}(S)$.

\vspace{2 mm}

\noindent
The set of all local transformations arising from such one-parameter families is obviously equal to $\Gamma(S)$. When \textit{E} is just an infinitesimal prolongation space, we can still define $\mathcal{L}(S)$ by using the tangency with respect to \textit{im S} and afterwards defining $\Gamma_{inf}(S)$ as being the pseudo-group generated (through composition and gluing) by the set of all the elements belonging to the local one-parameter families verifying $\varphi_{t_0}=Id~$ together with the above condition (b). We can further simplify (re-parametrize) by taking $t_0=0$ and clearly $\Gamma_{inf}(S)\subset \Gamma(S)$ when \textit{E} is also a finite prolongation space, $\mathcal{L}(S)$ being its infinitesimal generator. The condition $\varphi_{t_0}=Id$ is an acceptable simplification (by composing, if necessary, with $\varphi_{t_0}^{-1}$) and enables us to initiate the one-parameter family with the automorphism \textit{Id} that can always be prolonged.

\vspace{4 mm}

\noindent
Let us now examine the \textit{defining equations} for $\Gamma(S)$ and $\mathcal{L}(S)$ hence examine the nature of the infinitesimal jets that are automorphisms of \textit{S} up to a certain order and this brings us to examine \textit{contact} properties. In general, being given a differentiable manifold \text{M} and two sub-manifolds $N,~N'$ of the same dimension \textit{p} and meeting at a point \textit{x}, we shall say, according to the general definition of higher order Grassmannians, that these two sub-manifold have a $k-th$ order contact at the point \textit{x} when there exists an invertible $k-$jet $\mu$ of \textit{N} onto $N'$ with source and target \textit{x} such that $j_k\iota (x)=j_k\iota '(x)\cdot\mu$, where $\iota$ and $\iota '$ are the inclusions of the sub-manifolds. The equivalence classes under this relation are called the $k-th$ order contact elements of dimension \textit{p} on the manifold \textit{M} and at the point \textit{x}. Let us now take a vector field $\xi$ defined in a neighborhood of \textit{x}. The composite $\xi\circ\iota$ is a local section of the vector bundle $TM|_N\longrightarrow N$ , hence its $k-$jet at the point \textit{x} is an element of $J_k(TM|_N)$ and, of course, $J_kTN\subset J_k(TM|_N)$. We shall say that $\xi$ is tangent of order \textit{k} to \textit{N} at \text{x} when  $j_k(\xi\circ\iota)(x)\in J_kTN~$. Furthermore, if $N$ and $N'$ have a contact of order $k+1$ at the point \textit{x} then $TN$ and $TN'$ have, along any point of $T_x N=T_x N'$ a contact of order \textit{k} when considered as sub-manifolds of $TM$: The contact relation is obtained by extending to tangent vectors the initial contact relation. In particular, we infer that the vector fields tangent of order \textit{k} to the sub-manifold $N$ at \textit{x} are the same as those tangent to $N'$. Conversely, the last condition (on the tangency of $TN$ with $TN'$) implies the $(k+1)-st$ order tangency of $N$ with $N'$ at the point \textit{x} (\cite{Molino1977}, Propos.3, pg.20). Observing that $j_k(\xi\circ\iota)(x)=j_k\xi(x)\cdot j_k\iota(x)$ , we conclude that the $k-th$ order tangency concept for vector fields is more properly a notion of tangency of the elements of $J_kTM$ with a sub-manifold $N$ or, in a more general form, with the elements belonging to the fibre bundle of $(k+1)-st$ order $p-$dimensional contact elements of $M$ ($p=dim~N$). We finally remark that the groupoid $\Pi_k M$ operates to the left on the fibre bundle of $k-th$ order contact elements. There is an alternative definition, due to Ehresmann, stating that a $k-th$ order contact element is a linear subspace in the space of $k-th$ order tangent vectors. Though apparently simpler, this definitions seems to lack visibility. It is easy to \textit{see} tangency but not so easy to \textit{see} partial derivatives.

\vspace{4 mm}

\noindent
We now translate the above definitions of contact order into a form more suitable to our context, return to the local automorphisms and consider, to start with, a finite prolongation space $~E\longrightarrow P~$ of order $\ell$ as well as a structure \textit{S} of species \textit{E}. A local diffeomorphism $\varphi$ of \textit{P} is called a $k-th$ order automorphism of \textit{S} at the point \textit{y} whenever the image sub-manifold $p\varphi(im~S)$ has a $k-th$ order contact with \textit{im S} at the point $S\circ\varphi(y)$. Since the contact only depends on $j_k p\varphi(S(y))$ and the prolongation space is of finite order $\ell$, we infer that this notion boils down to the following: An element $Z\in\Pi_{\ell+k} P$ with source \textit{y} and target \textit{y'} is a $k-th$ order automorphism of \textit{S} when $pZ$, $k-$jet of source $S(y)$ corresponding to \textit{Z} ($(\ell+k)-$jets prolong to $k-$jets), transforms the $k-th$ order contact element of \textit{S} at $S(y)$ onto the $k-th$ order contact element of \textit{S} at $S(y')$. However, the contact element $pZ(S)$ is represented by the image of the local section $p\varphi\circ S\circ\varphi^{-1}$ where $\varphi$ is a representative of \textit{Z} i.e., $Z=j_{\ell+k}\varphi~$. Moreover, the images of two local sections $\sigma$ and $\sigma'$ of \textit{E} define the same $k-th$ order contact element at a common point $z\in E$ if and only if $j_k\sigma(\pi z)=j_k\sigma'(\pi z)$ and therefore $pZ(S)$ and \textit{S} define the same $k-th$ order contact element at $S(y')\in E$ if and only if $j_k(p\varphi\circ S\circ\varphi^{-1})(y')=j_k S(y')$ that is to say, when

\begin{equation*}
(pZ)\cdot j_kS(y)\cdot Z_k^{-1}=j_kS(y')~,\hspace{5 mm} Z_k=\rho_{k,\ell+k}Z .
\end{equation*}

\vspace{4 mm}

\newtheorem{highest}[LemmaCounter]{Lemma}
\begin{highest}
A jet $Z\in\Pi_{\ell+k}P$ is a $k-th$ order automorphism of S if and only if $Z(j_kS(\alpha Z))=j_kS(\beta Z)$ by means of the left action of $\Pi_{\ell+k}P$ on $J_kE$ (notations: $Z(~)=p_k\varphi(~)$, $\varphi$ represents Z and $p_k=p_k\circ p$). The $(\ell+k)-th$ order equation defining $\Gamma(S)$ is provided by the closed sub-groupoid (P being replaced by $\alpha(S)$ whenever S is not global)

\begin{equation}
\mathcal{R}_{\ell+k}(S)=\{Z\in\Pi_{\ell+k}P~|~Z(j_kS(\alpha Z))=j_kS(\beta Z)\}.
\end{equation}

\end{highest}

\vspace{4 mm}

\noindent
This equation is not \textit{complete} in general \textit{i.e.}, $\Gamma(S)$ might be strictly contained in $Sol~\mathcal{R}_{\ell+k}
(S)$ and, furthermore, the inclusion $\mathcal{R}_{\ell+k}(S)\supset J_{\ell+k}\Gamma(S)$ might also be strict.

\vspace{4 mm}

\noindent
We next consider the non-linear differential operator

\begin{equation*}
\mathcal{D}(S):~\varphi\in\underline{Di\!f\!f~P}~\longmapsto ~\varphi^{-1}(S)=(p\varphi^{-1})\circ S\circ\varphi\in\underline{E}
\end{equation*}

\vspace{2 mm}

\noindent
that is of order $\ell$ since it is obtained by the composition of $j_\ell$ and the fibered morphism

\begin{equation*}
\Phi(S):~Y\in\Pi_\ell P~\longmapsto~Y^{-1}(S(\beta Y))\in E
\end{equation*}

\vspace{2 mm}

\noindent
where $Y(~)$ is the action of the groupoid $\Pi_\ell P$ on \textit{E}. The $k-th$ order prolongation of $\mathcal{D}(S)$ is the operator of order $\ell+k$

\begin{equation*}
p_k\mathcal{D}(S):~\varphi\in\underline{Di\!f\!f~P}~\longmapsto ~j_k(p\varphi^{-1}\circ S\circ\varphi)=p_k\varphi^{-1}\circ j_kS\circ\varphi=
\end{equation*}
\begin{equation*}
\hspace{30 mm}=\varphi^{-1}(j_kS)\in\underline{J_kE}
\end{equation*}

\vspace{2 mm}

\noindent
and is obtained, with the aid of $j_{\ell+k}$, by means of the fibered morphism

\begin{equation*}
p_k\Phi(S):~Z\in\Pi_{\ell+k}P~\longmapsto~Z^{-1}(j_kS(\beta Z)\in J_kE
\end{equation*}

\vspace{2 mm}

\noindent
where $Z(~)$ is the action of $\Pi_{\ell+k}P$ on $J_kE$. The following short non-linear sequences

\begin{equation*}
1\longrightarrow\mathcal{R}_{\ell+k}(S)\longrightarrow\Pi_{\ell+k}P\xrightarrow{p_k\Phi(S)}(J_kE,j_kS)
\end{equation*}

\vspace{2 mm}

\noindent
and

\begin{equation*}
1\longrightarrow\underline{\Gamma(S)}\longrightarrow\underline{Di\!f\!f P}\xrightarrow{p_k\mathcal{D}(S)}(\underline{J_kE},j_kS)
\end{equation*}

\vspace{2 mm}

\noindent
are therefore exact in the set theoretical sense. When $im~p_k\Phi(S)$ is a sub-manifold of $J_kE$ and

\begin{equation*}
p_k\Phi(S):~\Pi_{\ell+k}P~\longrightarrow~im~p_k\Phi(S)
\end{equation*}

\vspace{2 mm}

\noindent
is a submersion then $\mathcal{R}_{\ell+k}(S)$ is a regularly embedded sub-manifold of $\Pi_{\ell+k}P$ and, for all \textit{h}, the equation $\mathcal{R}_{\ell+k+h}(S)$ is the prolongation in the usual sense i.e., the sub-set $~\Pi_{\ell+k+h}P\cap J_h\mathcal{R}_{\ell+k}(S)$ of $\mathcal{R}_{\ell+k}(S)~$. Taking into account the Proposition 2.1 in \cite{Gold1969}, we are tempted to replace the above two hypotheses by the unique assumption:

\begin{equation*}
p_k\Phi(S):~\Pi_{\ell+k}P~\longrightarrow~J_kE
\end{equation*}

\vspace{2 mm}

\noindent
is locally of constant rank in a neighborhood of each point belonging to $\mathcal{R}_{\ell+k}(S)$. Unfortunately (or perhaps fortunately) the above Proposition is inexact. If, with the notations of the above citation, X is reduced to a point, this Proposition would imply the remarkable statement: \textit{Théorème. Toute immersion est un plongement}. We shall return later to this matter and show that, in the specific case of the equations $\mathcal{R}_{\ell+k}(S)$, the first hypothesis can in fact be replaced by weaker conditions.

\vspace{4 mm}

\noindent
Let us now examine the infinitesimal aspects. A local vector field $\xi$ defined on \textit{P} is said to be a $k-th$ order infinitesimal automorphism of the structure \textit{S} at the point \textit{y} when the prolonged vector field $p\xi$ is tangent up to order \textit{k} to the sub-manifold $im~S$ at the point $z=S(y)$. Since $\pi\circ S=Id$, this condition can be replaced by

\begin{equation*}
j_k(p\xi\circ S)(y)=j_k(TS\circ\xi)(y),
\end{equation*}

\vspace{2 mm}

\noindent
this last condition measuring the "distancing" of $p\xi$ from $S_*\xi$ along $im~S$ in the vicinity of the point $z=S(y)$. The prolongation space \textit{E} being of finite order $\ell$, this last condition translates by the following: 
An element $Y\in J_{\ell+k}TP$ of source \textit{y} is a $k-th$ order infinitesimal automorphism of \textit{S} whenever 
\textit{pY}, $k-$jet of source $S(y)$ corresponding to \textit{Y}, is tangent of order \textit{k} to \text{im S} at the point
$S(y)$. We next observe that both $p\xi\circ S$ and $TS\circ\xi$ are local sections of the \textit{tangent Lie fibration} $TE\longrightarrow P$ composite of the natural projection $TE\longrightarrow E$ with $\pi$ and, consequently, the \textit{
reduced form} of the holonomic prolongation (\cite{Kumpera1971}, Théorème 12.1) shows that the vector $(p_k\xi)_{j_kS(y)}$ ($p_k=p_k\circ p$) only depends on $j_k(p\xi\circ S)(y)$ and the second line of the diagram (22.29) in the above citation, $\S$ 22, restricted to the fixed section $j_kS$ (the section being fixed, we are not forced any more to ascend to $J_{k+1}E$), shows that the jets $j_k(TS\circ\xi)(y)$ are precisely those, among the elements of $\tilde{J}_kTE$, for which the vector associated by prolongation is tangent to $im~j_kS$ at the point $j_kS(y)$. Summing up, we obtain the following

\newtheorem{lowest}[LemmaCounter]{Lemma}
\begin{lowest}
The jet $Y\in J_{\ell+k}TP$ is a $k-th$ order infinitesimal automorphism of S if and only if the vector $p_kY$ corresponding to it, at the point $j_kS(\alpha Y)$, is tangent to $im~j_kS$ ($p_kY=(p_k\xi)_{j_kS(\alpha Y)}$ where $\xi$ is a representative of \textit{Y}). Moreover, the equation of order $\ell+k$ of the Lie algebroid $\mathcal{L}(S)$ is given by

\begin{equation}
\textbf{R}_{\ell+k}(S)=\{Y\in J_{\ell+k}TP~|~p_kY\in (j_kS)_*T_{\alpha Y}P\}.
\end{equation}
\end{lowest}

\vspace{4 mm}

\noindent
Since the map $\lambda_k$ is a differentiable vector bundle morphism and the associated infinitesimal action $\Gamma(\lambda_k)$ is also a morphism of Lie algebra pre-sheafs, it follows that the linear equation \textbf{R}$_{\ell+k}(S)$ is a Lie equation (\cite{Kumpera1972}, \cite{Malgrange1972}). Moreover, $\mathcal{L}(S)=Sol~$\textbf{R}$_{\ell+k}(S)$ though, in general, this equation is not \textit{complete}. We next consider the linear differential operator
\begin{equation*}
\textbf{D}(S):~\xi\in\underline{TP}~\longmapsto~\theta(\xi)S=\frac{d}{dt} \varphi_t^{-1}(S)|_{t=0}\in\underline{VE~|~im~S}
\end{equation*}

\vspace{2 mm}

\noindent
where $VE$ is the vector sub-bundle of $TE$ composed by the $\pi-$vertical vectors and $\varphi_t$ is a local one-parameter family \textit{e.g.}, the local group, such that $\frac{d}{dt}\varphi_t|_{t=0}=\xi$. The restriction $VE~|~im S$ can be considered a bundle with base space \textit{P} since it can be identified to $S^{-1}VE$ and, according to the general definitions, $\bf{D}(S)$ is the linearization of $\mathcal{D}(S)$ along the section \textit{Id} of $\underline{Di\!f\!f~\!P}$. It is a linear operator of order $\ell$ whose associated linear morphism

\begin{equation*}
 \Psi(S):~J_\ell TP\longrightarrow VE~|~im~S
\end{equation*}

\vspace{2 mm}

\noindent
is defined as follows: Let $Y=j_\ell\xi(y),~(\varphi_t)$ a local one-parameter group associated to $\xi$ and $Y_t=j_\ell\varphi_t(y)$. Then
\begin{equation*}
\Psi(S)(Y)=\frac{d}{dt}[Y_t^{-1}(Sy_t)]_{t=0},
\end{equation*}

\vspace{2 mm}

\noindent
where $y_t=\beta Y_t$ and where, for simplicity, we write $Sy_t$ though meaning $S(y_t)$. However, $Y_t^{-1}(Sy_t)=p\varphi_t^{-1}(Sy_t)=p\varphi_{-t}(Sy_t)$ and, consequently, $\Psi(S)(X)=-(p\xi)_{S(y)}+S_*(\xi_y)$. We infer that

\begin{equation}
\Psi(S)(Y)=-\textbf{v}(p\xi_{S(y)}),
\end{equation}

\vspace{2 mm}

\noindent
where \textbf{v}( ) is the vertical component of a vector following the direct sum decomposition $T_{S(y)}E=V_{S(y)}E~\oplus~S_*(T_yP)$ along \textit{im S}.

\vspace{2 mm}

\noindent
Let us now return to the general definitions of \cite{Kumpera1975}, I, $\S$17 and, in particular, to the last sequence, on pg.341, exhibiting the linear morphism  $\tau_k(T\mathcal{D})$ that defines the differential operator $T\mathcal{D}$. In the present case (and adapting slightly the notations), it concerns the vertical linearization along the section $Id$ for the specific fibration $P=M~\times~M~\longrightarrow~M~$, $(x,x')~\mapsto~x$ , where $J_\ell P$ is replaced by $J_\ell M=J_\ell(M,M)$ (each section of $P$ identifies with a map $M~\longrightarrow~M$), $\tilde{J}_\ell VP~|~J_\ell M=J_\ell M~\times_M~J_\ell TM$ (fibre product with respect to the projection $\beta$) and where $p_\ell:~J_\ell M~\times_M~J_\ell TM~\longrightarrow~VJ_\ell M$ simply becomes the canonical identification. The morphism

\begin{equation*}
T\tau_\ell(\mathcal{D})=T\Phi(S):~VJ_\ell M~|~im~Id~\longrightarrow~VE~|~im~S
\end{equation*}

\vspace{2 mm}

\noindent
becomes, by means of the canonical identification, equal to $\Psi(S)$. The $k-th$ order prolongation of $\textbf{D}(S)$ is then the linear operator of order $\ell+k$

\begin{equation*}
p_k\textbf{D}(S):~\xi\in\underline{TM}~\longmapsto~j_k[\theta(\xi)S]\in\underline{J_k(VE~|~im~S)},
\end{equation*}

\vspace{2 mm}

\noindent
that we can, by means of the isomorphism $p_k:~\tilde{J}_k VE~\longrightarrow~VJ_kE$ (\cite{Kumpera1975}, Propos.12.2), replace by

\begin{equation*}
\underline{p_k}\circ p_k\textbf{D}(S):~\underline{TM}~\longrightarrow~\underline{VJ_k E~|~im~j_kS}
\end{equation*}

\vspace{2 mm}

\noindent
that, in turn, is nothing else but the linearized, along the identity section, of $p_k\mathcal{D}(S)$ (\cite{Kumpera1975}, pg.342). Furthermore, the linear morphism associated to $\underline{p_k}\circ p_k\textbf{D}(S)$ is on the one hand equal to $p_k\circ p_k\Psi(S)$ and on the other, for being a linearization and due to the canonical identification, equal to
\begin{equation}
T\tau_{\ell+k}(p_k\mathcal{D})=Tp_k\Phi(S):~VJ_{\ell+k}M~|~im~Id~\longrightarrow~VJ_k E~|~im~S.
\end{equation}

\vspace{2 mm}

\noindent
An entirely similar calculation, where we shall replace 

\begin{equation*}
\varphi_t^{-1}(S)=\mathcal{D}(S)\varphi_t\hspace{4 mm}by\hspace{4 mm}\varphi_t^{-1}(j_kS)=p_k\mathcal{D}(S)\varphi_t,
\end{equation*}

\vspace{2 mm}

\noindent
will show that

\begin{equation}
p_k\circ p_k\Psi(S)(X)=-\textbf{v}(p_k\xi_{j_kS(x,x')})
\end{equation}

\vspace{2 mm}

\noindent
where $X\in J_{\ell+k}TM$ and $v$ is the vertical component of a vector in the direct sum decomposition $T_ZJ_kE=V_ZJ_kE~
\oplus~(j_kS)_*(T_xM),~Z=j_kS(x,x')$ and this implies the exactness of the sequence

\begin{equation}
0~\longrightarrow~\textbf{R}_{\ell+k}(S)~\longrightarrow~J_{\ell+k}TM~\xrightarrow{p_k\Psi(S)}~J_k(VE~|~im~S)
\end{equation}

\vspace{2 mm}

\noindent
for all $k\geq 0~$, hence $\textbf{R}_{\ell+k+h}(S)$ is the prolongation of order \textit{h} of the equation $\textbf{R}_{\ell+k}(S)$ whenever $p_k\Psi(S)$ is locally of constant rank as a linear mapping on each fibre (\textit{i.e.}, when $\textbf{R}_{\ell+k}(S)$ is a locally trivial vector bundle). Moreover, the sequence

\begin{equation}
0~\longrightarrow~\underline{\mathcal{L}(S)}~\longrightarrow~\underline{TM}~\xrightarrow{p_k\textbf{D}(S)}~\underline{J_k(VE~|~im~S)}
\end{equation}

\vspace{2 mm}

\noindent
is also exact for all $k\geq 0$.

\section{Local and infinitesimal equivalences}
In this section we study the consequences of the previous hypotheses and conditions in view of showing that the \textit{micro-differentiable structures} introduced by Pradines (\cite{Pradines1966}) and hereafter considered, do generate the desired global structures. By the canonical identification $\Pi_{\ell+k}P~\times_P~J_{\ell+k}TP~\longrightarrow~V\Pi_{\ell+k}P$ resulting from the prolongation, by the target, of local vector fields defined on \textit{P}, each sub-space $\textbf{R}_{\ell+k}(S)_y$ determines, at every point $X\in\Pi_{\ell+k}P$ with $\beta(X)=y$, a sub-space $(\Delta_{\ell+k})_X\subset V_X\Pi_{\ell+k}P$ of the same dimension as that of $\textbf{R}_{\ell+k}(S)$ and thus defines an $\alpha-$vertical distribution on $\Pi_{\ell+k}P$ (field of contact elements) whose point-wise dimension is locally constant if and only if the fibre bundle $\textbf{R}_{\ell+k}(S)$ is locally trivial. This last condition being satisfied, the distribution $\Delta_{\ell+k}$ is integrable since $\textbf{R}_{\ell+k}(S)$ is a Lie equation.  

\vspace{2 mm}

\newtheorem{integer}[PropositionCounter]{Proposition}
\begin{integer}
Let S be a structure of species E and finite order $\ell$. Then $ker~Tp_k\Phi(S)=\Delta_{\ell+k}(S)$ and, furthermore, the following properties are equivalent:

\vspace{2 mm}

\hspace{2 mm}i) $p_k\Psi(S)$ is locally of constant rank,

\vspace{2 mm}

\hspace{1 mm}ii) $p_k\Phi(S)$ is locally of constant rank along the units of $\Pi_{\ell+k}P$,

\vspace{2 mm}

iii) $p_k\Phi(S)$ is locally of constant rank, 

\vspace{2 mm}

\hspace{1 mm}iv) $\Delta_{\ell+k}$ is locally of constant dimension.

\vspace{2 mm}

\noindent
These equivalent conditions being verified, $\textbf{R}_{\ell+k}(S)$ is a union of integral leaves of the distribution $\Delta_{\ell+k}$.

\end{integer}

\noindent
$\bf{Proof.}$ We first observe that $\Delta_{\ell+k}$ is invariant by all the right translations of the groupoid $\Pi_{\ell+k}P$, such a translation mapping $\alpha-$fibres onto $\alpha-$fibres and preserving the targets. On the other hand, the morphism $p_k\Phi(S)$ is a \textit{differential co-variant} relative to the right action of the groupoid $\Pi_{\ell+k}P$ on itself \textit{i.e.}, the following formula holds:
\begin{equation}
p_k\Phi(S)(X\cdot Y)=Y^{-1}p_k\Phi(S)(X).
\end{equation}

\vspace{4 mm}

\noindent
We next observe that $ker~Tp_k\Phi(S)~|~im~Id=\Delta_{\ell+k}~|~im~Id$ since $p_k\circ p_k\Psi(S)$ identifies, by means of the canonical identification, to $Tp_k\Phi(S)~|~im~Id=\Delta_{\ell+k}~|~im~Id$ (\text{cf.}(14)) and that the kernel of $Tp_k\Phi(S):~T\Pi_{\ell+k}P~\longrightarrow~TJ_kE$ is equal to that of the restriction $Tp_k\Phi(S):~V\Pi_{\ell+k}P~\longrightarrow~VJ_kE~$. The invariance of $\Delta_{\ell+k}$ and the co-variance of $p_k\Phi(S)$ entails the equality everywhere and we infer the equivalence of the four stated properties. If $\mathcal{F}$ is a leaf of $\Delta_{\ell+k}$, then $Tp_k\Phi(S)~|~T\mathcal{F}=0$ and consequently $p_k\Phi(S)(\mathcal{F})$ reduces to a point. In particular, when $X\in\mathcal{F}$ and $p_k\Phi(X)=j_kS(x),~x=\alpha(X)$, then $p_k\Phi(S)(\mathcal{F})=j_kS(x)$ where after $\mathcal{F}\subset\mathcal{R}_{\ell+k}(S)$ and this achieves the proof.

\vspace{4 mm}

\noindent
Let us now observe that, for any fibration morphism $\lambda$,
\begin{equation*}
P~\xrightarrow{~\lambda~}~P
\end{equation*}
\begin{equation*}
\downarrow\hspace{12 mm}\downarrow
\end{equation*}
\begin{equation*}
M\xrightarrow{~Id~}M
\end{equation*}

\vspace{2 mm}

\noindent
the following two conditions are equivalent:

\vspace{2 mm}

\textit{a) $\lambda$ is locally of constant rank,}

\vspace{2 mm}

\textit{b) $\lambda$ is locally of vertical constant rank i.e., the rank of the restrictions of $\lambda$ to the fibres is locally constant with respect to the topology of P (and not only of the fibres),}

\vspace{2 mm}

\noindent
since we always have the relation $rank_y\lambda=dim~M+(vertical~rank)_y\lambda~$. This amounts to say that the rank of $T\lambda:~TP~\longrightarrow~TP'$ differs, at each point, from the rank of $T\lambda:~VP~\longrightarrow~VP'$ by the integer \textit{dim M} and, in particular, that $ker~T\lambda=ker~T\lambda~|~VP$. We also observe that the above properties still hold when we replace \textit{Id} by any diffeomorphism $\varphi$ of \textit{M} (and even by a diffeomorphism $M\longrightarrow M'$). In particular, when \textit{P} and \textit{P'} are fibre bundles and $\lambda$ is a morphism of such bundles, then the vertical rank (rank of $T\lambda$ restricted to the tangent space of a fibre) is equal to the rank of the  restriction of $\lambda$ to the fibres.

\vspace{2 mm}

\noindent
On the other hand, and without any regularity hypotheses pending upon $p_k\Phi(S)$ or $p_k\Psi(S)$, we remark that the isotropy

\begin{equation}
(\mathcal{R}^0_{\ell+k}S)_y=\{X\in\mathcal{R}_{\ell+k}(S)~|~\alpha(X)=\beta(X)=y\}
\end{equation}

\vspace{2 mm}

\noindent
of $\mathcal{R}_{\ell+k}(S)$ at the point \textit{y} is always a closed Lie sub-group of the total isotropy $(\Pi^0_{\ell+k}M)_y$ since it is given by the "closed" conditions $X(j_kS(y))=j_kS(y)$. Its Lie algebra identifies canonically with the linear isotropy

\begin{equation}
(\textbf{R}^0_{\ell+k}S)_y=\{X\in\textbf{R}_{\ell+k}(S)~|~\alpha(X)=y~,~\beta(X)=0\}
\end{equation}

\vspace{2 mm}

\noindent
of $\textbf{R}_{\ell+k}(S)$ that on its turn is determined by the condition $p_kX=0~$. Furthermore, these two remarks show that the finite and infinitesimal $k-th$ order isotropies of \textit{S} at the point \textit{y} are entirely determined by the jet $j_kS(y)$. In particular, the jet of order $\ell$ only depends on the point $S(y)\in E$. Let us denote by $\mathcal{R}_{\ell+k}(S)_y$ the fibre, with respect to $\alpha$ and above the point \textit{y}, of $\mathcal{R}_{\ell+k}(S)$ and by $\mathcal{B}_{\ell+k}(S)_y$ its projection $\beta(\mathcal{R}_{\ell+k}(S)_y)$.

\vspace{2 mm}

\newtheorem{isotropy}[CorollaryCounter]{Corollary}
\begin{isotropy}
If \textbf{R}$_{\ell+k}(S)$ is a locally trivial vector bundle then, at any point y, the fibre $\mathcal{R}_{\ell+k}(S)_y$ is a closed and regularly embedded sub-manifold of $(\Pi_{\ell+k}P)_y$ whose connected components are integral leaves of $\Delta_{\ell+k}$. Furthermore, $\mathcal{B}_{\ell+k}(S)_y$ is a closed and regularly embedded sub-manifold of P and the triple $(\mathcal{R}_{\ell+k}(S)_y~,~\beta~,~\mathcal{B}_{\ell+k}(S)_y)$ is a locally trivial sub-fibre bundle of $(\Pi_{\ell+k}P)_y~|~\mathcal{B}_{\ell+k}(S)_y$ with structure group equal to $\mathcal{R}^0_{\ell+k}(S)_y$ (it being understood that $\alpha(S)=P$). 
\end{isotropy}

\noindent
$\bf{Proof.}$ It is clear that $\mathcal{R}_{\ell+k}(S)_y$ is closed being the inverse image of the point $j_kS(y)$ relative to the map $p_k\Phi(S)_y:~(\Pi_{\ell+k}P)_y~\longrightarrow~(J_kE)_y$. Inasmuch, it is also clear that $\mathcal{R}_{\ell+k}(S)_y$ is a principal space of the isotropy group $\mathcal{R}^0_{\ell+k}(S)_y$, the orbits being the inverse images, by $\beta$, of the points of $\mathcal{B}_{\ell+k}(S)_y~$. The hypothesis on $\textbf{R}_{\ell+k}(S)$ entails, in virtue of the previous proposition, that $p_k\Phi(S)_y:~(\Pi_{\ell+k}P)_y~\longrightarrow~(J_kE)_y$ has a locally constant rank and consequently $\mathcal{R}_{\ell+k}(S)_y$ is a regularly embedded sub-manifold. Since $ker~Tp_k\Phi(S)=\Delta_{\ell+k}(S)$, we see at once that $T\mathcal{R}_{\ell+k}(S)_y=\Delta_{\ell+k}(S)~|~\mathcal{R}_{\ell+k}(S)_y$ and thus infer that the leaves of $\Delta_{\ell+k}(S)$ contained in $\mathcal{R}_{\ell+k}(S)_y$ are open sets hence, due to the connexity of the leaves, are necessarily the connected components. We denote by $\Xi$ the distribution defined on \textit{P} by $\Xi_y=\beta(\textbf{R}_{\ell+k}(S)_y)$. The right invariance of the distribution $\Delta_{\ell+k}$ or, still better, the definition itself of this distribution shows immediately that $\beta_*(\Delta_{\ell+k})_X=\Xi_{\beta (X)}$. This distribution $\Xi$ is not, in general, of locally constant dimension but is generated by a family of vector fields that is involutive and locally of finite type. In fact, every section of $\textbf{R}_{\ell+k}(S)$ determines, by projection, a vector field that is a section of $\Xi$ and, since $\textbf{R}_{\ell+k}(S)$ is a Lie equation, the bracket of two sections projects onto the bracket of their images. Moreover, since  $\textbf{R}_{\ell+k}(S)$ is locally trivial, the pre-sheaf of local sections of this fibre bundle is locally of finite type (in fact, locally free) and the local finiteness property of $\Xi$ follows. It is also easy to verify, by using again the right translations of $\Pi_{\ell+k}P$, that the integral leaves of $\Xi$ are precisely the projections of the integral leaves of $\Delta_{\ell+k}$ (\textit{cf.} \cite{Hermann1962} and \cite{Turiel1976}, Chap.I, $\S~2$).
Hence, we thus conclude that $\beta(\mathcal{R}_{\ell+k}(S)_y)=\mathcal{B}_{\ell+k}(S)_y$ is a union of integral leaves of $\Xi$ and this union is discrete: For every $z\in\mathcal{B}_{\ell+k}(S)_y$, there exists an open neighborhood \textit{U} of \textit{z} in \textit{P} such that $\mathcal{B}_{\ell+k}(S)_y\cap U$ reduces to the intersection of \textit{U} with a unique integral leaf of $\Xi$. To see this, it suffices to note firstly that the leaves of $\Delta_{\ell+k}$ passing by the points of the same $\beta-$fibre of $\mathcal{R}_{\ell+k}(S)_y$ project all on the same leave of $\Xi$ and secondly, recalling that $\mathcal{R}_{\ell+k}(S)_y$ is a regularly embedded sub-manifold whose connected components are integral leaves of $\Delta_{\ell+k}~$, we shall take an open neighborhood $\mathcal{U}$, in $\Pi_{\ell+k}(S)_y~$, of a point \textit{X} contained in a $\beta-$fibre of $\mathcal{R}_{\ell+k}(S)_y$ in such a way that $\mathcal{U}\cap\mathcal{R}_{\ell+k}(S)_y$ just contains the points of a single leaf of $\Delta_{\ell+k}$. By right translation with the elements of the isotropy $\mathcal{R}^0_{\ell+k}(S)_y$, the same situation reproduces itself, with the translated open set, at every other point of the $\beta-$fibre and consequently $U=\beta(\mathcal{U})$ responds to the required property. Shrinking, if necessary, the open set $\mathcal{U}$ and recalling the regularity of the embedding of $\mathcal{R}_{\ell+k}(S)_y$, we can show further that $\mathcal{B}_{\ell+k}(S)_y\cap U$ is a \textit{slice} (in a coordinate system) and, consequently, that $\mathcal{B}_{\ell+k}(S)_y$ is a regularly embedded sub-manifold of \textit{P}. Finally, a saturation argument of $\mathcal{R}_{\ell+k}(S)_y$ by the action of the total isotropy  $(\Pi^0_{\ell+k}P)_y$ and the fact that $\mathcal{R}_{\ell+k}(S)_y$ is closed, shows that $\mathcal{B}_{\ell+k}(S)_y$ is closed in \textit{P}. We thus see that $\beta:\mathcal{R}_{\ell+k}(S)_y~\longrightarrow~\mathcal{B}_{\ell+k}(S)_y$ is a \textit{surmersion} (surjective submersion) and, consequently, a (sub-)principal fibre bundle of $(\Pi_{\ell+k}P)_y~|~\mathcal{B}_{\ell+k}(S)_y$ that is locally trivial since it admits local sections and the proof is therefore complete. We shall nevertheless continue by showing, further, that the dimension of the eventually non-connected sub-manifold $\mathcal{B}_{\ell+k}(S)_y$ is constant on each connected component of $\alpha(S)$. In fact, since the dimension of $\mathcal{R}_{\ell+k}(S)$ is constant on each connected component $\mathcal{O}$ of $\alpha(S)$, we infer that the dimension of $\Delta_{\ell+k}(S)$ is constant on $\beta^{-1}\mathcal{O}$ and, consequently, that of $\mathcal{R}_{\ell+k}(S)_y$ is also constant above $\mathcal{O}$. On the other hand, the tangent space to each $\beta-$fibre of $\mathcal{R}_{\ell+k}(S)_y$ is isomorphic to the Lie algebra $\textbf{R}^0_{\ell+k}(S)_y$ of the structural group $\mathcal{R}^0_{\ell+k}(S)_y$ hence the dimension of this tangent space is constant. Lastly, since $\Xi=T_y\mathcal{B}_{\ell+k}(S)_y$ is isomorphic to the quotient of $(\Delta_{\ell+k})_Y~,~Y\in\mathcal{R}_{\ell+k}(S)_y~,~\beta(Y)=y~$, modulo the tangent space to the $\beta-$fibre at the point \textit{Y}, the result follows. To terminate, we provide an alternative proof of the above statements in view of the  geometrical mechanisms that it will turn apparent and that will be of relevance in the sequel. For this, we go back to the definition of the elements of $\textbf{R}_{\ell+h}(S)$ as being the jets $j_{\ell+h}\xi(y)$ of local vector fields $\xi$ whose prolongations $p\xi$ are $h-th$ order tangent to the section \textit{S} at the point $S(y)$. However, we can see that any $Y\in\mathcal{R}_{\ell+k}(S)$ transforms $\mathcal{R}_{\ell+h}(S)_{\alpha Y}$ into $\mathcal{R}_{\ell+h}(S)_{\beta Y}$, $h<k$ , since the $k-$jet $pY$ with source $S(\alpha Y)$ and associated, by prolongation, to \textit{Y} transforms the $k-th$ order contact element of $im~S$ at the point $S(\alpha Y)$ into the corresponding contact element at the point $S(\beta Y)$. This property, however, fails when $h=k$ since $pY$ only operates on  $(k-1)-$jets of vector fields on \textit{E} and the invariance does not subsist any longer, not even for the projected sub-spaces $\beta(\textbf{R}_{\ell+h}(S)_{\alpha Y})$ and $\beta(\textbf{R}_{\ell+h}(S)_{\beta Y})$. However, if we observe in general, \textit{N} being an arbitrary manifold, that $\Pi_rN$ operates on $J^0_rTN~$, then it will become apparent (\textit{cf.} Lemmas 3 and 4) that \textit{Y} transforms  $\textbf{R}^0_{\ell+k}(S)_{\alpha Y}$ onto $\textbf{R}^0_{\ell+k}(S)_{\beta Y}$ since the jet $j_k(p\xi)[S(\alpha Y)]$ associated to $j_{\ell+k}\xi(\alpha Y)\in\textbf{R}^0_{\ell+k}(S)_{\alpha Y}$ belongs to $(J^0_kTE)_{S(\alpha Y)}$ and, consequently, $pY$ transforms the jet $j_k(p\xi)[S(\alpha Y)]$, $k-th$ order tangent to $im~S$ at the point $S(\alpha Y)$, into a $k-$jet of vector field tangent, up to order \textit{k}, to $im~S$ at the point $S(\beta Y)$, this later $k-$jet being precisely the one that corresponds, via prolongation, to the transformed jet $Y(j_{\ell+k}\xi(\alpha Y))\in\textbf{R}^0_{\ell+k}(S)_{\beta Y}~$. We thus infer that the fibres of $\textbf{R}^0_{\ell+k}(S)$ at the points $\alpha Y$ and $\beta Y$ are isomorphic where after the isomorphy of all the fibres along any orbit $\mathcal{B}_{\ell+k}(S)_y$ of $\textbf{R}_{\ell+k}(S)$ in \textit{P}. Since $\Xi_y=\textbf{R}_{\ell+k}(S)_y~/~\textbf{R}^0_{\ell+k}(S)_y$, we infer that $\Xi_y$ has constant dimension along $\mathcal{B}_{\ell+k}(S)_y$ if and only if it is inasmuch for $\textbf{R}_{\ell+k}(S)_y~$. In particular, this implies the constancy of dimensions for $\Xi_y$ on the intersection of $\mathcal{B}_{\ell+k}(S)_y$ with a connected component of $\alpha S$. The sub-manifold $\mathcal{B}_{\ell+k}(S)_y$ has therefore a constant dimension in each connected component of $\alpha(S)$.

\vspace{4 mm}

\noindent
The method of proof suggests a weakening of the regularity hypotheses imposed on $\textbf{R}_{\ell+k}(S)$. Accordingly, it would suffice to assume that $\textbf{R}_{\ell+k}(S)$ is locally of finite type \textit{i.e.}, that it be locally generated by a finite number of sections (differentiable sections of $J_{\ell+k}TP$ taking values in $\textbf{R}_{\ell+k}(S)$). The bracket properties of $\underline{J_{\ell+k}TP}$ would then imply that $\textbf{R}_{\ell+k}(S)$ is a Lie equation, the set of its local sections being closed under the bracket. However, in the specific case of the equation $\textbf{R}_{\ell+k}(S)$, this generalization is only illusory. In fact, the local finiteness of generators would imply lower semi-continuity for  $dim~\textbf{R}_{\ell+k}(S)_y$ whereas the definition of this equation as the kernel of $p_k\Psi(S)$ would imply the upper semi-continuity. In definite, the local finiteness assumption is entirely equivalent to regularity.

\section{Micro-differentiable structures and globalization}
In this section we look for the hypotheses enabling us to endow the groupoid $\mathcal{R}_{\ell+k}(S)$ or eventually its $\alpha-$connected component with a differentiable structure compatible with its algebraic structure. To begin, we assume that the vector bundle $\textbf{R}_{\ell+k}(S)$ is locally trivial and already possesses all the regularity properties indicated in the last Proposition as well as in its Corollary. We denote by $\mathcal{R}_{\ell+k}(S)_0$ the union of all the integral leaves of $\Delta_{\ell+k}$ issued from the units of $\Pi_{\ell+k}P~$. A standard connectivity argument shows that $\mathcal{R}_{\ell+k}(S)_0$ is a sub-groupoid of $\mathcal{R}_{\ell+k}(S)$ that we shall call its $\alpha-$connected component of the units space (assumed to be connected otherwise we argue on each connected component). Every $\alpha-$fibre of $\mathcal{R}_{\ell+k}(S)_0$ is in fact the connected component of a unit in the corresponding $\alpha-$fibre of $\mathcal{R}_{\ell+k}(S)$.

\vspace{4 mm}

\noindent
According to the general definitions (\cite{Kumpera1971}, \cite{Pradines1967}), $\mathcal{R}_{\ell+k}(S)_0$ is the sub-groupoid of $\Pi_{\ell+k}P$ generated by the Lie algebroid $\underline{\textbf{R}_{\ell+k}(S)}\subset\underline{J_{\ell+k}TP}~$. Contrary to what has been said and written in the last century (\cite{Rodrigues1962}, main theorem), it is well known that Lie's Second Theorem is inexact even for transitive Pseudo-groups \textit{i.e.}, for transitive sub-groupoids of the general groupoid $\Pi_{\ell+k}P~$. In other terms, given a Lie sub-algebroid $\mathcal{A}$ of $\underline{J_kTP}$ and denoting by $\mathcal{G}$ the (algebraic) sub-groupoid of $\Pi_{\ell+k}P$ it generates (\textit{e.g.}, by integrating the $\alpha-$fibres distribution), it is not always possible to endow this sub-groupoid with a differentiable structure (of sub-manifold) in such a way that its Lie algebroid can be identified with the given one. The main obstruction rests in the non-vanishing \textit{holonomy} for the integral leaves of $\Delta_{\ell+k}~$, the distribution generated by the right translations applied to $\mathcal{A}~$, these leaves being precisely the $\alpha-$fibres of the desired sub-groupoid. Nonetheless, in our present situation, this holonomy fortunately vanishes since $\mathcal{R}_{\ell+k}(S)$ is the kernel of a differential operator (or, more precisely, the kernel of its defining morphism). We shall of course proceed in the most standard way by first endowing an open neighborhood of the units with the differentiable structure practically "imposed" by the algebroid $\underline{\textbf{R}_{\ell+k}(S)}$ and, thereafter, propagate this local differentiable structure to the $\alpha-$connected component $\mathcal{R}_{\ell+k}(S)_0~$. In order to further propagate this differentiable structure to the whole of $\mathcal{R}_{\ell+k}(S)$, we shall be forced to add an additional hypothesis. Let us also observe that we are undertaking this painstaking homework since, to our knowledge, this construction has never been fully elucidated before.

\vspace{4 mm}

\noindent
Since $\Delta_{\ell+k}$ is regular and integrable, we take, for each unit \textbf{e} in $\mathcal{R}_{\ell+k}(S)$, a \textit{foliating chart}  ($\mathcal{U}_\alpha,\phi_\alpha)$ (\cite{Chevalley1946}, p.69) in such a way that $\textbf{e}\in\mathcal{U}_\alpha$ and that the leaves contained in $\mathcal{U}_\alpha$ are \textit{slices} with respect to the coordinate system $\phi_\alpha~$, where a slice means a sub-manifold diffeomorphic to an open $\textbf{p}-$cube in $\bf{R}^p$, $\textbf{p}=dim~\Delta_{\ell+k}$ (\textit{cf.} the aforementioned reference). Let us denote by $\mathcal{V}_\alpha$ the union of all the slices contained in $\mathcal{U}_\alpha$ and that contain units of $\mathcal{R}_{\ell+k}(S)$. Then $\mathcal{V}_\alpha$ is a closed and regularly embedded sub-manifold of $\mathcal{U}_\alpha~$. We set $\mathcal{U}=\bigcup~\mathcal{U}_\alpha$ and $\mathcal{V}=\bigcup~\mathcal{V}_\alpha~$. Then $\mathcal{V}_\alpha\cap\mathcal{V}_\beta$ is an open sub-set of both $\mathcal{V}_\alpha$ and $\mathcal{V}_\beta~$. Consequently, $\mathcal{V}$ is a closed and regularly embedded sub-manifold of $\mathcal{U}$ and, of course, $\mathcal{V}\subset\mathcal{R}_{\ell+k}(S)_0~$. We next observe that the projection $\alpha:\mathcal{V}\longrightarrow\alpha(S)$ is a surmersion and that the $\alpha-$fibres contained in $\mathcal{V}$ are sub-manifolds of locally constant dimension with respect to the variation of \textit{y} in $\alpha(S)$. Moreover, the germ of this sub-manifold, along the units, is uniquely determined \textit{i.e.}, does not depend on the initial choice of the foliating charts and, by construction, $\textbf{R}_{\ell+k}(S)$ can be identified with the $\alpha-$vertical tangent bundle of $\mathcal{V}$ along the units.

\vspace{4 mm}

\noindent
Due to its geometrical interest, we shall construct the same germ of sub-manifold by a procedure relying on the local constancy of the rank of $p_k\Phi(S)$. A unit $\textbf{e}\in\mathcal{R}_{\ell+k}(S)$ being fixed, there exists, due to the constancy of the rank, a neighborhood $\mathcal{U}$ of \textbf{e} in $\Pi_{\ell+k}P$ such that $\mathcal{W}=p_k\Phi(S)(\mathcal{U})$ is a regularly embedded sub-manifold of $J_kE$ and that $p_k\Phi(S):\mathcal{U}~\longrightarrow~\mathcal{W}$ is a surmersion. Furthermore, we can assume that $\overline{\mathcal{U}}=\alpha(\mathcal{U})=\alpha(\mathcal{U}\cap\mathcal{I})$, where $\mathcal{I}$ represents the sub-manifold of units of $\Pi_{\ell+k}P$ and, thereafter, the intersection $\mathcal{W}\cap im~j_kS=j_kS(\overline{\mathcal{U}})$ is a regularly embedded sub-manifold of $\mathcal{W}~$. The inverse image of $j_kS(\overline{\mathcal{U}})$ by the map $p_k\Phi(S)~|_{\mathcal{U}}$ is equal to $\mathcal{R}_{\ell+k}(S)\cap\mathcal{U}$ and consequently, this inverse image is a closed and regularly embedded sub-manifold $\mathcal{V}$ of $\mathcal{U}$. Since $ker~Tp_k\Phi(S)=\Delta_{\ell+k}~$, we infer that the $\alpha-$fibres of $\mathcal{V}$ are integral manifolds of maximum dimension of $\Delta_{\ell+k}~$, these $\alpha-$fibres being precisely the fibres of $p_k\Phi(S):\mathcal{U}~\longrightarrow~\mathcal{W}$ above the points of $j_kS(\overline{\mathcal{U}})$. By shrinking, if necessary, the open set $\mathcal{U}$, we can be brought to the case where these fibres are slices and thus infer that any leave of $\Delta_{\ell+k}~$, issued from a point in $\mathcal{U}\cap\mathcal{I}~$, intercepts the open set $\mathcal{U}$ along a unique slice. Furthermore, the open set $\mathcal{U}$ is a foliating chart defined in a neighborhood of the unit \textbf{e} and verifies the conditions stated previously. By taking the union of all these open sets $\mathcal{U}_\alpha$ corresponding to the various units, we obtain an open neighborhood $\tilde{\mathcal{U}}$ of the units of $\mathcal{R}_{\ell+k}(S)$ in $\Pi_{\ell+k}P$ such that $\mathcal{V}=\tilde{\mathcal{U}}\cap\mathcal{R}_{\ell+k}(S)=\tilde{\mathcal{U}}\cap\mathcal{R}_{\ell+k}(S)_0$ is the union of the slices contained in the open sets $\mathcal{U}_\alpha$ that contain the units. Besides, $\mathcal{V}$ is a neighborhood of the units in  $\mathcal{R}_{\ell+k}(S)$ as well as in  $\mathcal{R}_{\ell+k}(S)_0~$. The germ, along these units, of the regularly embedded sub-manifold $\mathcal{V}$ is unique and will be called, together with Pradines (\cite{Pradines1966}), the \text{micro-differentiable structure} of $\mathcal{R}_{\ell+k}(S)$.

\vspace{4 mm}

\noindent
We now show that this micro-differentiable structure can be extended to a global structure defined on $\mathcal{R}_{\ell+k}(S)_0~$. Indeed, if $Y\in\mathcal{R}_{\ell+k}(S)_0$ and if $\textbf{e}=\alpha(Y)$, the continuation Theorem (\cite{Palais1957}, pg.10) enables us to define a differentiable mapping $\mu:\overline{\mathcal{U}}~\longrightarrow~\Pi_{\ell+k}P$ where $\overline{\mathcal{U}}$ is an open neighborhood of the unit \textbf{e} in the sub-manifold $\mathcal{I}$ composed by the units, $\mu(\textbf{e'}),~\textbf{e'}\in\overline{\mathcal{U}}~$, belongs to the leaf of $\Delta_{\ell+k}$ that contains \textbf{e'} ($\alpha-$fibre of $\mathcal{R}_{\ell+k}(S)_0$) and $\mu(\textbf{e})=Y$. However, if we identify \textit{P} with the units manifold $\mathcal{I}~$, $\mu$ becomes a differentiable section of the bundle $\alpha:\Pi_{\ell+k}P~\longrightarrow~P$ defined on $\overline{\mathcal{U}}$ with values in $\mathcal{R}_{\ell+k}(S)_0$ that assumes the value \textit{Y} at the point \textbf{e}. Due to the local rank constancy of $p_k\Phi(S)$, there exists a neighborhood $\mathcal{U}$ of \textit{Y} in $\Pi_{\ell+k}P$ such that $\mathcal{W}=p_k\Phi(S)(\mathcal{U})$ is a regularly embedded sub-manifold of $J_kE$ and that $p_k\Phi(S):\mathcal{U}~\longrightarrow~\mathcal{W}$ is a surmersion. We can further suppose that $\overline{\mathcal{U}}=\alpha(\mathcal{U})=\alpha(\mathcal{U}\cap im~\mu)$ which implies that $\mathcal{W}\cap im~j_kS=j_kS(\overline{\mathcal{U}})$ is a regularly embedded sub-manifold of $\mathcal{W}$. It then follow as previously and shrinking, if necessary, the open set $\mathcal{U}$, that $\mathcal{R}_{\ell+k}(S)\cap\mathcal{U}=\mathcal{R}_{\ell+k}(S)_0\cap\mathcal{U}$ is a closed and regularly embedded sub-manifold of $\mathcal{U}$, inverse image by the map $p_k\Phi(S)|_{\mathcal{U}}$ of $j_kS(\overline{\mathcal{U}})$, and the fibres of the fibration $p_k\Phi(S)|_{\mathcal{U}}$ coincide, above $j_kS(\overline{\mathcal{U}})$, with the $\alpha-$fibres of $\mathcal{R}_{\ell+k}(S)\cap\mathcal{U}~\longrightarrow~\overline{\mathcal{U}}~$. It follows therefore that $\mathcal{R}_{\ell+k}(S)_0$ is endowed with the differentiable structure of a regularly embedded sub-manifold of $\Pi_{\ell+k}P$ compatible with its groupoid structure. Moreover, there exists an open set $\mathcal{U}$ in $\Pi_{\ell+k}P$ such that $\mathcal{R}_{\ell+k}(S)_0$= $\mathcal{R}_{\ell+k}(S)\cap\mathcal{U}$ or, in other terms, that $\mathcal{R}_{\ell+k}(S)_0$ is locally closed. The open set $\mathcal{U}$ can be chosen saturated with respect to the leaves of $\Delta_{\ell+k}$ since $\mathcal{R}_{\ell+k}(S)$ as well as $\mathcal{R}_{\ell+k}(S)_0$ are already saturated. The constructions show clearly that $\textbf{R}_{\ell+k}(S)$ identifies with the $\alpha-$vertical tangent bundle of $\mathcal{R}_{\ell+k}(S)_0$ along the units and that, consequently, the sheaf $\underline{\textbf{R}_{\ell+k}(S)}$ is the Lie algebroid associated to the differentiable sub-groupoid $\mathcal{R}_{\ell+k}(S)_0$ (\cite{Kumpera1971},\cite{Pradines1967}).

\newtheorem{different}[PropositionCounter]{Proposition}
\begin{different}
Let S be a structure of species E and finite order $\ell$ verifying the equivalent conditions of the Proposition $1.$ The equation $\mathcal{R}_{\ell+k}(S)_0$ is then an $\alpha-$connected and regularly embedded Lie sub-groupoid of $\Pi_{\ell+k}P$ whose associated Lie algebroid is equal to $\underline{\textbf{R}_{\ell+k}(S)}$. There exists an open set $\mathcal{U}$ of $\Pi_{\ell+k}P~$, saturated by the integral foliation of $\Delta_{\ell+k}~$, such that $\mathcal{R}_{\ell+k}(S)_0$ is closed in $\mathcal{U}$ and, furthermore, the sequence

\begin{equation}
1~\longrightarrow~\mathcal{R}_{\ell+k}(S)_0~\longrightarrow~\mathcal{U}~\xrightarrow{p_k\Phi(S)}~(J_kE,j_kS)
\end{equation}

\vspace{2 mm}

\noindent
is exact. Conversely, when $\mathcal{R}_{\ell+k}(S)_0$ admits a structure of Lie sub-groupoid of $\Pi_{\ell+k}P$ with associated Lie algebroid equal to $\underline{\textbf{R}_{\ell+k}(S)}$, then the equivalent conditions of the Proposition $1$ are satisfied and the differentiable structure of $\mathcal{R}_{\ell+k}(S)_0$ is regularly embedded.  
\end{different}

\vspace{4 mm}

\noindent
This proposition states precisely that $\mathcal{R}_{\ell+k}(S)_0$ is the non-linear Lie equation generated by the linear Lie equation $\textbf{R}_{\ell+k}(S)$ (\textit{cf.} \cite{Malgrange1970},\cite{Malgrange1972}) when $p_k\Psi$ is locally of constant rank. The converse statement of the above proposition is quite simple to prove and will be omitted. Nevertheless, it veils behind the curtains the following rather important fact: We first observe that any local section of $\textbf{R}_{\ell+k}(S)$ determines, by right translations, a right invariant vector field on $\Pi_{\ell+k}P$ that is contained in $\Delta_{\ell+k}$ and conversely. Let $\mathcal{P}$ denote the pseudo-group of local transformations operating on $\Pi_{\ell+k}P$ and generated by the flows of the above mentioned right invariant vector field. Then, each element of $\mathcal{P}$ leaves invariant the sub-groupoid $\mathcal{R}_{\ell+k}(S)_0$ though it needs not preserve the larger sub-groupoid $\mathcal{R}_{\ell+k}(S)$, this being the main obstacle towards the possibility of simply extending or prolonging the differentiable structure of the former to the later, as was done in the micro-differentiable situation.
So, let us now take care of $\mathcal{R}_{\ell+k}(S)$.

\vspace{4 mm}

\noindent
In the previous attempt to provide $\mathcal{R}_{\ell+k}(S)_0$ with a differentiable structure, there is a key point that still remains unexplored: The open sets $\mathcal{U}~$, neighborhood of the unit \textbf{e} , can be chosen in such a way that any leaf of $\Delta_{\ell+k}$ issued from a point $\textbf{e'}\in\mathcal{U}\cap\mathcal{I}$ cuts the open set $\mathcal{U}$ along a unique slice namely, the one containing \textbf{e'}. We show in fact that $\mathcal{U}$ can be chosen in such a way that any leaf of $\Delta_{\ell+k}$ meets $\mathcal{U}$ at most along a single slice (no holonomy). Since the operations of the groupoid $\Pi_{\ell+k}M$ are continuous, we can take open neighborhoods $\mathcal{V}$ and $\mathcal{W}$ of \textbf{e} such that $\mathcal{V}\subset\mathcal{U}~$, $\mathcal{W}\subset\mathcal{U}~$, $\mathcal{V}\cdot\mathcal{V}^{-1}\subset\mathcal{U}$ and $\mathcal{W}\cdot\mathcal{V}\subset\mathcal{U}~$, the operations being performed on all composable pairs. We can further assume that both neighborhoods are the domains of foliating charts, each fibre being a slice. Let us now take a leaf $\mathcal{F}$ and assume that it meets  $\mathcal{V}$ along two slices $\mathcal{S}_1$ and $\mathcal{S}_2~$. 
If we take a point $X\in\mathcal{S}_1$ then $\mathcal{F}_{\beta(X)}=\mathcal{F}\cdot X^{-1}$ is the leaf of $\Delta_{\ell+k}$ passing by the unit $\beta(X)\in\mathcal{U}\cap\mathcal{I}$ and consequently $\mathcal{s}_1\cdot X^{-1}$ as well as $\mathcal{S}_2\cdot X^{-1}$ are included in the slice of $\mathcal{U}$ hence also in that of $\mathcal{W}$ and containing the unit $\beta(X)$. Applying the inverse operation $\cdot X~$, we infer that $\mathcal{S}_1$ and $\mathcal{S}_2$ are both contained in a same slice of $\mathcal{U}$ and consequently in a same slice of $\mathcal{V}$ since the later is just a foliated chart restriction of $\mathcal{U}~$. More generally, we show that the same property continues to hold at each point of $\Pi_{\ell+k}P$, it being understood that \textit{P} is to be replaced by $\alpha(S)$ when \textit{S} is just a local section. 
In fact, let $X=j_{\ell+k}\varphi(y)$ be an arbitrary point and let us consider the \textit{flow} $j_{\ell+k}\varphi~$. By right translations, provided by the flow elements, we establish a diffeomorphism $\tau:\beta^{-1}(\beta\varphi)~\longrightarrow~\beta^{-1}(\alpha\varphi)~$ that is compatible with $\Delta_{\ell+k}$ and consequently the leaves are transformed in leaves, the foliating charts in foliating charts and the slices in slices. Furthermore, the desired property is verified for the foliating open set $\tau(\mathcal{U})$ as soon as it is verified for $\mathcal{U}~$. In sum, for every $X\in\Pi_{\ell+k}P~$, there exists a foliating chart $\mathcal{U}$ of $\Delta_{\ell+k}$ that is a neighborhood of \textit{X} and for which any leaf of $\Delta_{\ell+k}$ meets at most along a single slice. However, these properties translate by saying that the integral foliation of $\Delta_{\ell+k}$ is \textit{simple} or, in other terms, (\cite{Palais1957}, pg.19) that there exists a differentiable structure, necessarily unique, on the quotient $\Pi_{\ell+k}P~/~\Delta_{\ell+k}$ of the general groupoid modulo the integral leaves of $\Delta_{\ell+k}$ and in such a way that the quotient map $\zeta:\Pi_{\ell+k}M~\longrightarrow~\Pi_{\ell+k}M~/~\Delta_{\ell+k}~$ is a surmersion. Moreover, since each leaf of $\Delta_{\ell+k}$ is contained in an $\alpha-$fibre, there is a canonical projection $\overline{\alpha}$ of the quotient space onto \textit{P} that commutes with $\alpha$. As previously, we shall also denote by $\mathcal{I}$ the identity section $y\in P~\longmapsto~j_{\ell+k}Id(y)\in\Pi_{\ell+k}P$ and set $\overline{\mathcal{I}}=\zeta\circ\mathcal{I}~$. Then, of course, $\overline{\mathcal{I}}$ is a differentiable section of $\overline{\alpha}~$, the restriction $\zeta:im~\mathcal{I}~\longrightarrow~im~\overline{\mathcal{I}}$ is a diffeomorphism, $im~\overline{\mathcal{I}}$ is a regularly embedded sub-manifold of the quotient $\Pi_{\ell+k}P~/~\Delta_{\ell+k}$ and $\mathcal{R}_{\ell+k}(S)_0=\zeta^{-1}(im~\overline{\mathcal{I}})$. Moreover, since $\mathcal{R}_{\ell+k}(S)_0$ is locally closed in a saturated open set $\mathcal{U}~$, the sub-manifold $im~\overline{\mathcal{I}}$ is locally closed in the open set $\zeta(\mathcal{U})$ (it should however be observed that $\Pi_{\ell+k}P~/~\Delta_{\ell+k}$ needs not be separated). Finally, since $\mathcal{R}_{\ell+k}(S)$ and $\mathcal{R}_{\ell+k}(S)_y$ are closed and saturated by the leaves and since $\mathcal{R}_{\ell+k}(S)_y$ is a regularly embedded sub-manifold of $(\Pi_{\ell+k}P)_y$ whose connected components are precisely the maximal integral leaves of $\Delta_{\ell+k}~$, we conclude that $\zeta(\mathcal{R}_{\ell+k}(S))$ is closed and that $\zeta(\mathcal{R}_{\ell+k}(S)_y)$ is closed and discrete (each point is isolated) in the fibre $(\Pi_{\ell+k}P~/~\Delta_{\ell+k})_y~$.

\vspace{4 mm}

\noindent
Let us next assume that $\mathcal{R}_{\ell+k}(S)$ is endowed with a differentiable structure that makes it become a Lie sub-groupoid of $\Pi_{\ell+k}P$ and whose associated Lie algebroid is equal to $\underline{\textbf{R}_{\ell+k}(S)}$. We show that, under these conditions, the differentiable structure of $\mathcal{R}_{\ell+k}(S)$ is regularly embedded in $\Pi_{\ell+k}P~$, $\zeta(\mathcal{R}_{\ell+k}(S))$ becomes a regularly embedded sub-manifold of $\Pi_{\ell+k}P~/~\Delta_{\ell+k}$ admitting the open subset $im~\overline{\mathcal{I}}$ and the restriction $\overline{\alpha}:\zeta(\mathcal{R}_
{\ell+k}(S))~\longrightarrow~\alpha(S)$ is \textit{étale}. 

\vspace{4 mm}

\noindent
In fact, since $\mathcal{R}_{\ell+k}(S)$ is a Lie groupoid, the projection

\begin{equation*}
\alpha:\mathcal{R}_{\ell+k}(S)~\longrightarrow~\alpha(S) 
\end{equation*}

\vspace{2 mm}

\noindent
is a surmersion and consequently, for any $X\in\mathcal{R}_{\ell+k}(S)$, there exists a differentiable local section $\mu$ of $\alpha$ taking its values in $\mathcal{R}_{\ell+k}(S)$ and passing through \textit{X}. Furthermore, there exists an open neighborhood $\mathcal{V}$ of \textit{X} in $\mathcal{R}_{\ell+k}(S)$ such that the fibres of $\alpha:\mathcal{V}~\longrightarrow~\alpha(\mathcal{V})$ are slices. By the hypothesis, $\underline{\textbf{R}_{\ell+k}(S)}$ is the Lie algebroid of $\mathcal{R}_{\ell+k}(S)$ hence the $\alpha-$vertical tangent bundle $V\mathcal{R}_{\ell+k}(S)$ is equal to $\Delta_{\ell+k}|\mathcal{R}_{\ell+k}(S)$. We infer that the slices of $\mathcal{V}$ are integral sub-manifolds of maximum dimension of $\Delta_{\ell+k}$ and, more generally that the $\alpha-$fibres of $\mathcal{R}_{\ell+k}(S)$ have, for connected components, the integral leaves of $\Delta_{\ell+k}~$. Let $\mathcal{F}$ be the leaf that meets $X=\mu(y)$, let us take a second $\alpha-$section $\nu$ of $\mathcal{R}_{\ell+k}(S)$ such that $\nu(y)\in\mathcal{F}$ and let us denote by $\tilde{\mathcal{V}}$ the saturated set of $\mathcal{V}$, in $\mathcal{R}_{\ell+k}(S)$, by the integral leaves of $\Delta_{\ell+k}~$. Furthermore, the continuation Theorem implies that $\tilde{\mathcal{V}}$ is an open subset of $\mathcal{R}_{\ell+k}(S)$ and, since $\nu(y)\in\tilde{\mathcal{V}}$, we infer that $\nu(z)$ and $\mu(z)$ belong to the same leaf of $\Delta_{\ell+k}$ as soon as \textit{z} is sufficiently close to \textit{y}. However, this implies that the two sections $\zeta\circ\mu$ and $\zeta\circ\nu$ of $\overline{\alpha}$ coincide in a neighborhood of \textit{y} hence enables to define a structure of differentiable sub-manifold on the image $\zeta(\mathcal{R}_{\ell+k}(S))$ in such a way that $\overline{\alpha}:\zeta(\mathcal{R}_{\ell+k}(S))~\longrightarrow~\alpha(S)$ becomes \textit{étale}. Let us finally show that this sub-manifold is regularly embedded. To do so, we recall that $\mathcal{R}_{\ell+k}(S)$ is defined as being the kernel of $p_k\Phi(S)$ and that this mapping is locally of constant rank. Therefore, we can find an open neighborhood $\mathcal{U}$ of \textit{X} in $(\Pi_{\ell+k}P)$ such that $\mathcal{W}=p_k\Phi(S)(\mathcal{U})$ is a regularly embedded sub-manifold of $J_kE$ and $p_k\Phi(S):\mathcal{U}~\longrightarrow~\mathcal{W}$ is a surmersion. We can further assume that $\overline{\mathcal
{U}}=\alpha(\mu)=\alpha(\mathcal{U})=\alpha(\mathcal{U}\cap im~\mu)$ and this implies that $\mathcal{W}\cap im~j_kS=j_kS(\overline{\mathcal{U}})$ is a regularly embedded sub-manifold of $\mathcal{W}$. We infer that $\mathcal{R}_{\ell+k}(S)\cap\mathcal{U}$ is a closed and regularly embedded sub-manifold of $\mathcal{U}$ , inverse image of $j_kS(\overline{\mathcal{U}})$ by the map $p_k\Phi(S)|_{\mathcal{U}}$, and whose $\alpha-$fibres coincide with the fibers of  $p_k\Phi(S)|_{\mathcal{U}}$ above the points of $j_kS(\overline{\mathcal{U}})$. Since, $ker~p_k\Phi(S)=\Delta_{\ell+k}$ , the $\alpha-$fibres of $\alpha:\mathcal{R}_{\ell+k}(S)\cap\mathcal{U}~\longrightarrow~\overline{\mathcal{U}}$ are integral sub-manifolds of maximal dimension of $\Delta_{\ell+k}$ and consequently, it is possible to choose the above open set $\mathcal{U}$ as well as the open neighborhood $\mathcal{V}$ of \textit{X} in $\mathcal{R}_{\ell+k}(S)$, considered in the beginning, in such a way that $\mathcal{V}=\mathcal{R}_{\ell+k}(S)\cap\mathcal{U}$ and, moreover, that the two differentiable structures, one induced by the given structure of $\mathcal{R}_{\ell+k}(S)$ and the other induced by that of $\mathcal{U}$ , be the same and shows consequently that the given structure on $\mathcal{R}_{\ell+k}(S)$ is regularly embedded. Saturating these open sets by the leaves of $\Delta_{\ell+k}$ , we obtain much in the same way the open and saturated sub-sets 

\begin{equation*}
\tilde{\mathcal{V}}=\mathcal{R}_{\ell+k}(S)\cap\tilde{\mathcal{U}}
\end{equation*}

\noindent
and

\begin{equation*}
\zeta(\mathcal{R}_{\ell+k}(S)\cap\mathcal{U})=\zeta(\mathcal{R}_{\ell+k}(S)\cap\tilde{\mathcal{U}})=\zeta(im~\mu\cap\tilde{\mathcal{U}})=im~(\zeta\circ\mu)\cap\zeta(\tilde{\mathcal{U}})
\end{equation*}

\vspace{2 mm}

\noindent
which implies that $\zeta(\mathcal{R}_{\ell+k}(S))$ is a regularly embedded sub-manifold of $\Pi_{\ell+k}P~/~\Delta_{\ell+k}~$.

\vspace{2 mm}

\noindent
Conversely, when $\zeta(\mathcal{R}_{\ell+k}(S))$ admits the structure of a regularly embedded sub-manifold, then $\mathcal{R}_{\ell+k}(S)$ admits the structure of a regularly embedded sub-manifold of $\Pi_{\ell+k}P$ that induces forcefully the regularly embedded structure of $\mathcal{R}_{\ell+k}(S)_ y~$. Consequently, the Lie algebroid associated to the Lie sub-groupoid $\mathcal{R}_{\ell+k}(S)$ is necessarily equal to $\underline{\textbf{R}_{\ell+k}(S)}$ and this implies again that the structure of $\zeta(\mathcal{R}_{\ell+k}(S))$ is \textit{étale} over $\alpha(S)$. It is further clear that $im~\overline{\mathcal{I}}$ is an open sub-manifold of $\zeta(\mathcal{R}_{\ell+k}(S))$ that however needs not be closed by virtue of the eventual non-separability of the quotient $\Pi_{\ell+k}P~/~\Delta_{\ell+k}$ as well as that of $\zeta(\mathcal{R}_{\ell+k}(S))$.

\vspace{2 mm}

\noindent
We show as well, using the same arguments as above, that if every $X\in\mathcal{R}_{\ell+k}(S)$, the latter without any previously assigned  structure, belongs to the image of a differentiable $\alpha-$section of $\Pi_{\ell+k}P$ taking its values in $\mathcal{R}_{\ell+k}(S)$, then $\mathcal{R}_{\ell+k}(S)$ is a regularly embedded sub-manifold of $\Pi_{\ell+k}P~$. Such sections are, by the way,  usually obtained by composing sections of $\overline{\alpha}:\zeta(\mathcal{R}_{\ell+k}(S))~\longrightarrow~\alpha(S)$ with sections of $\zeta~$. In sum, we proved the following:

\vspace{2 mm}

\newtheorem{equivalent}[PropositionCounter]{Proposition}
\begin{equivalent}
Let S be a structure of species E and finite order $\ell$ satisfying the equivalent requirements of the Proposition $1$. Under these conditions, the following properties are also equivalent:

\vspace{2 mm}

\hspace{2 mm}i) Every element of $\mathcal{R}_{\ell+k}(S)$ belongs to the image of some differentiable section of $\alpha:\Pi_{\ell+k}P~\longrightarrow~P$ taking its values in $\mathcal{R}_{\ell+k}(S)$.

\vspace{2 mm}

\hspace{1 mm}ii) The closed subset $\zeta(\mathcal{R}_{\ell+k}(S))\subset\Pi_{\ell+k}P~/~\Delta_{\ell+k}$ is a regularly embedded sub-manifold \textit{étale} over $\alpha(S)$. 

\vspace{2 mm}

iii) $\mathcal{R}_{\ell+k}(S)$ is a Lie sub-groupoid of $\Pi_{\ell+k}P$ with associated Lie algebroid equal to $\underline{\textbf{R}_{\ell+k}(S)}$.

\vspace{2 mm}

\hspace{1 mm}iv) $\mathcal{R}_{\ell+k}(S)$ is a regularly embedded Lie sub-groupoid of $\Pi_{\ell+k}P~$.
\end{equivalent}

\vspace{4 mm}

It is desirable, at this point, to examine more carefully the embedding of $\mathcal{R}_{\ell+k}(S)$. With this in mind, let us assume that the structure \textit{S} satisfies the equivalent conditions of the last Proposition whereby $\mathcal{R}_{\ell+k}(S)$ becomes a regularly embedded 
Lie sub-groupoid of $\Pi_{\ell+k}P$ and let us denote by $\Pi_{\ell+k}P~/~\mathcal{R}_{\ell+k}(S)$ the set of right sided classes of $\Pi_{\ell+k}P$ modulo $\mathcal{R}_{\ell+k}(S)$, that is to say, the quotient modulo the relation: $X\sim Y$ if and only if $X\cdot Y^{-1}\in\mathcal{R}_{\ell+k}(S)$. These classes are no other than the orbits, by the left action (via the target), of $\mathcal{R}_{\ell+k}(S)$ on $\Pi_{\ell+k}P$ and consequently we see promptly that the right action (via the source) of $\Pi_{\ell+k}P$ on itself factors to an action (to the right) of $\Pi_{\ell+k}P$ on the above quotient $\Pi_{\ell+k}P~/~\mathcal{R}_{\ell+k}(S)$. We have already seen that each $\alpha-$fibre of $\mathcal{R}_{\ell+k}(S)$ is a closed and regularly embedded sub-manifold of an $\alpha-$fibre of $\Pi_{\ell+k}P$ and its connected components are the integral leaves of $\Delta_{\ell+k}~$. Let us now take a unit \textbf{e}. Since there exists a saturated open set $\tilde{\mathcal{U}}$ such that $\tilde{\mathcal{U}}\cap\mathcal{R}_{\ell+k}(S)=\tilde{\mathcal{U}}\cap\mathcal{R}_{\ell+k}(S)_0$ , we infer, in view of the previous results, that there exists a foliating chart $\mathcal{U}$ for $\Delta_{\ell+k}$ , neighborhood of \textbf{e}, such that each $\alpha-$fibre of $\mathcal{R}_{\ell+k}(S)$ issued from a point $\textbf{e'}\in\mathcal{U}\cap\mathcal{I}$ meets the open set $\mathcal{U}$ along a unique slice namely, the slice containing \textbf{e'}. Let us next observe that the right action of $\Pi_{\ell+k}P$ permutes the trajectories (orbits) of $\mathcal{R}_{\ell+k}(S)$. Arguments entirely analogous to those used previously for the integral foliation of $\Delta_{\ell+k}$ will show that, for any $X\in\Pi_{\ell+k}P$, there exists a foliating chart $\mathcal{V}$ of $\Delta_{\ell+k}$ , neighborhood of \textit{X} , such that an arbitrary trajectory of $\mathcal{R}_{\ell+k}(S)$ will meet the open set $\mathcal{V}$ in at most one slice. We find ourselves within conditions entirely analogous to those found in the Theorem 8, pg.19 of \cite{Palais1957}. In this theorem, the differentiable structure of the quotient (\textit{i.e.}, the appropriate changes of charts) is guaranteed by the transport Theorem (\cite{Palais1957}, pg.10) which however does not apply in the present case in view of the (eventual) non-connectivity of the $\alpha-$fibres of $\mathcal{R}_{\ell+k}(S)$. Nevertheless, the transport can be replaced by the following argument: Let $X\in\mathcal{R}_{\ell+k}(S)$, $\textbf{e}=\alpha(X)$ the corresponding source, $\mu$ a differentiable section of $\alpha:\mathcal{R}_{\ell+k}(S)~\longrightarrow~\alpha(S)$ assuming the value $\mu(\textbf{e})=X$ (it is essentially here that intervenes the property (\textit{i}) of the last Proposition) and $\mathcal{U}$ , respectively $\mathcal{U}'$, the domain of a foliating chart of $\Delta_{\ell+k}~$, neighborhood of \textbf{e}, respectively \textit{X}, whose intersection with any orbit of $\mathcal{R}_{\ell+k}(S)$ reduces at most to a single slice and that verifies moreover the condition $\mu(\mathcal{U}\cap\mathcal{I})\subset\mathcal{U}'$. Let $\mathcal{W}$ be a sub-manifold transverse to the slices of $\mathcal{U}$ such that $\mathcal{W}\supset\mathcal{U}\cap\mathcal{I}$ and let us set

\begin{equation*}
\mathcal{V}=\{Y\in\mathcal{W}~|~\beta(Y)\in\mathcal{U}\cap\mathcal{I}\}.
\end{equation*}

\vspace{4 mm}

\noindent
Clearly, $\mathcal{V}$ remains a transversal sub-manifold and we define the mapping $\Sigma:\mathcal{V}~\longrightarrow~\Pi_{\ell+k}P~,~\Sigma(Y)=\mu(\beta Y)\cdot Y$, and we see readily that the following properties hold:

\vspace{4 mm}

1) $\Sigma(Y)$ is contained in the orbit of \textit{Y},

\vspace{3 mm}

2) $\Sigma~|~\mathcal{U}\cap\mathcal{I}=\mu~$,

\vspace{3 mm}

3) $\Sigma(\mathcal{V})\subset\mathcal{U}'$, in shrinking if necessary the sub-manifold $\mathcal{V}$.

\vspace{4 mm}

\noindent
The condition (1) implies that $\Sigma$ is injective. Furthermore, since the right action of $\Pi_{\ell+k}P$ on itself is effective and the map

\begin{equation*}
\mu:~\mathcal{U}\cap\mathcal{I}~\longrightarrow~\mu(\mathcal{U}\cap\mathcal{I})
\end{equation*}

\vspace{2 mm}

\noindent
is a diffeomorphism of regularly embedded sub-manifolds, we infer that $\Sigma$ has injective rank (on the tangent level) and consequently $\Sigma:\mathcal{V}~\longrightarrow~\Sigma(\mathcal{V})$ is a diffeomorphism of sub-manifolds respecting the orbits (property (1)) and $\Sigma(\mathcal{V}$ is transverse to the slices of $\mathcal{U}'$. We can therefore extend the conclusions of Palais' Theorem to the space of the orbits of $\mathcal{R}_{\ell+k}(S)$ since the two charts in the quotient space originating from $\mathcal{U}$ and $\mathcal{U}'$ do as well originate from the transverse sub-manifolds $\mathcal{V}$ and $\Sigma(\mathcal{V})$ and therefore are compatible. The quotient set $\Pi_{\ell+k}P~/~\mathcal{R}_{\ell+k}(S)$ admits therefore a (necessarily unique) differentiable manifold structure for which the quotient map is a surmersion. We next remark that $\Pi_{\ell+k}P~/~\Delta_{\ell+k}=\Pi_{\ell+k}P~/~\mathcal{R}_{\ell+k}(S)_0$ and the diagram below is commutative, 

\begin{equation*}
\Pi_{\ell+k}P~/~\mathcal{R}_{\ell+k}(S)_0~\xleftarrow{\hspace{2 mm}\zeta\hspace{2 mm}}~\Pi_{\ell+k}P~\xrightarrow{~\zeta'~}~\Pi_{\ell+k}P~/~\mathcal{R}_{\ell+k}(S)
\end{equation*}
\begin{equation*}
\downarrow\hspace{30 mm}\downarrow\hspace{30 mm}\downarrow
\end{equation*}
\begin{equation*}
P\hspace{9 mm}\xleftarrow{~Id~}\hspace{9 mm}P\hspace{9 mm}\xrightarrow{~Id~}\hspace{9 mm}P
\end{equation*}

\vspace{4 mm}

\noindent
the arrow $\Upsilon=\zeta'\circ\zeta^{-1}$ being surjective and \textit{étale}. Furthermore, $\mathcal{R}_{\ell+k}(S)=(\zeta')^{-1}(\zeta'(\mathcal{I}))$ (inverse image) and $\zeta(\mathcal{R}_{\ell+k}(S))=\Upsilon^{-1}(\zeta'(\mathcal{I}))$. The arrows $\zeta,~\zeta'$ and $\Upsilon$ are differential co-variants with respect to the right action (by the source) of $\Pi_{\ell+k}P$ on the three spaces. The quotient manifold $\Pi_{\ell+k}P~/~\mathcal{R}_{\ell+k}(S)$ can be obtained from $\Pi_{\ell+k}P~/~\mathcal{R}_{\ell+k}(S)_0$ by identifying the points on each orbit that are deducible one from the other by the discrete action of $\mathcal{R}_{\ell+k}(S)$ and this identification is globally compatible when the equivalent properties of the above Proposition are verified. We next observe that the co-variance of $p_k\Phi(S)$ enables the factorisation of this morphism, modulo the action of $\mathcal{R}_{\ell+k}(S)$, and consequently the diagram that follows is commutative, the factored differential co-variant $p'_k\Phi(S)$ becoming an injective immersion (\textit{inmersion? Kkkk}).

\begin{equation*}
\Pi_{\ell+k}P~\xrightarrow{p_k\Phi(S)}~J_kE
\end{equation*}
\begin{equation}
\hspace{7 mm}\zeta'\downarrow\hspace{16 mm}\nearrow~p'_k\Phi(S)
\end{equation}
\begin{equation*}
\Pi_{\ell+k}P~/~\mathcal{R}_{\ell+k}(S)\hspace{5 mm}
\end{equation*}

\vspace{4 mm}

\noindent
We infer that $im~p_k\Phi(S)$ is a sub-manifold (not always regularly embedded) of $J_kE$ canonically isomorphic to the quotient $\Pi_{\ell+k}P~/~\mathcal{R}_{\ell+k}(S)$, the map $\alpha:im~p_k\Phi(S)~\longrightarrow~\alpha(S)$ is a fibration and

\begin{equation*}
p_k\Phi(S):~\Pi_{\ell+k}P~\longrightarrow~im~p_k\Phi(S)
\end{equation*}

\vspace{2 mm}

\noindent
is a surmersive morphism of fibrations having for basis $\alpha(S)$ (\textit{P} being replaced by $\alpha(S)$ when \textit{S} is not global). The covariance of $p_k\Phi(S)$ finally shows that $(im~p_k\Phi(S),\alpha,\alpha(S),p_k)$ is a prolongation space of order $\ell+k$ and that the map

\begin{equation*}
p'_k\Phi(S):(\Pi_{\ell+k}P~/~\mathcal{R}_{\ell+k}(S),\overline{\alpha},\alpha(S),p_k^S)\longrightarrow(im~p_k\Phi(S),\alpha,\alpha(S),p_k)
\end{equation*}

\vspace{2 mm}

\noindent
is an isomorphism of prolongation spaces where the first term is given the quotient prolongation space structure, modulo the right action of $\mathcal{R}_{\ell+k}(S)$, of the canonical structure of $\Pi_{\ell+k}P$ resulting from the standard prolongation operation by the source (\cite{Kumpera1975}, $\S$ 16, part (b)). 

\vspace{2 mm}

\newtheorem{quotient}[TheoremCounter]{Theorem}
\begin{quotient}
Let S be a structure of species E and of finite order $\ell$ such that $p_k\Psi(S)$ is locally of constant rank. Then the following conditions are equivalent:

\vspace{2 mm}

\hspace{2 mm}i) $\mathcal{R}_{\ell+k}(S)$ is a Lie sub-groupoid of $\Pi_{\ell+k}P$ whose associated Lie algebroid is equal to $\underline{\textbf{R}_{\ell+k}(S)}$.

\vspace{2 mm}

\hspace{1 mm}ii) $\mathcal{R}_{\ell+k}(S)$ is a regularly embedded (and closed) Lie sub-groupoid of $\Pi_{\ell+k}P~$.

\vspace{2 mm}

iii) There exists a differentiable structure on $\Pi_{\ell+k}P~/~\mathcal{R}_{\ell+k}(S)$ such that the quotient map is a submersion.

\vspace{2 mm}

\hspace{1 mm}iv) The image of $p_k\Phi(S)$ admits a sub-manifold structure of $J_kE$ such that the map $p_k\Phi(S):\Pi_{\ell+k}P~\longrightarrow~im~p_k\Phi(S)$ is a submersion.

\vspace{2 mm}

\hspace{2 mm}v) Every element of $\mathcal{R}_{\ell+k}(S)$ belongs to the image of a local differentiable section of $\alpha:\Pi_{\ell+k}P~\longrightarrow~P$ taking values in $\mathcal{R}_{\ell+k}(S)$.

\vspace{2 mm}

\noindent
These equivalent conditions being verified, the quotient differential covariant

\begin{equation*}
p'_k\Phi(S):(\Pi_{\ell+k}P~/~\mathcal{R}_{\ell+k}(S),\overline{\alpha},\alpha(S),p_k^S)\longrightarrow(im~p_k\Phi(S),\alpha,\alpha(S),p_k)
\end{equation*}

\vspace{2 mm}

\noindent
is an isomorphism of prolongation spaces.
\end{quotient}

\vspace{4 mm}

\newtheorem{sub-quotient}[CorollaryCounter]{Corollary}
\begin{sub-quotient}
Let S be a structure of species E and of finite order $\ell$ such that $p_k\Psi(S)$ is locally of constant rank and let us assume further that the equation $\textbf{R}_{\ell+k}(S)$ is transitive $(\beta(\textbf{R}_{\ell+k}(S))=TP)$. Under these conditions, the equivalent properties of the previous Theorem are always satisfied namely, $\mathcal{R}_{\ell+k}(S)$ is a closed and regularly embedded Lie sub-groupoid of $\Pi_{\ell+k}P$. Moreover, $\mathcal{R}_{\ell+k}(S)$ as well as $\mathcal{R}_{\ell+k}(S)_0$ are locally trivial Lie sub-groupoids and $\mathcal{R}_{\ell+k}(S)_0$ 
is closed in $\Pi_{\ell+k}P$. 
\end{sub-quotient}

\vspace{2 mm}

\noindent
$\bf{Proof.}$ The transitivity of $\textbf{R}_{\ell+k}(S)$ implies that the restriction of $\beta$ to each $\alpha-$fibre of $\mathcal{R}_{\ell+k}(S)$ is a submersion (that will be surjective whenever $\mathcal{R}_{\ell+k}(S)$ becomes transitive). Let $X\in\mathcal{R}_{\ell+k}(S)$, $y=\alpha(X)$, and let us take a section $\sigma$ of $\beta:\mathcal{R}_{\ell+k}(S)_x~\longrightarrow~P$ defined on an open neighborhood $\mathcal{U}$ of \textit{y} such that $\sigma(y)=\textbf{e}$ (the unit associated to \textit{y}). The map $\tau:\mathcal{U}~\longrightarrow~\Pi_{\ell+k}P~$, 
$\tau(y)=X\cdot\sigma(y)^{-1}$, is an $\alpha-$section taking its values in $\mathcal{R}_{\ell+k}(S)$ and such that $\tau(y)=X~$, thus retrieving the property (\textit{v}) of the last Theorem. The submersivity of $\beta$ on each $\alpha-$fibre of $\mathcal{R}_{\ell+k}(S)$ implies the submersivity of $\alpha\times\beta:\mathcal{R}_{\ell+k}(S)~\longrightarrow~\alpha(S)\times\alpha(S)$ whereupon the possibility (\cite{Kumpera1971}) in defining local trivialisations of $\mathcal{R}_{\ell+k}(S)$ and $\mathcal{R}_{\ell+k}(S)_0$ with the help of local sections of $\alpha\times\beta$. We also remark, \textit{en passant}, that each domain of a local trivialisation admits a regularly embedded differentiable structure inherited from the isotropy, at a given point, of $\mathcal{R}_{\ell+k}(S)$ (resp. $\mathcal{R}_{\ell+k}(S)_0$) , this isotropy being a closed Lie sub-group hence regularly embedded (resp. regularly embedded hence also closed). The family of all such trivialisations is compatible and defines thereafter the regularly embedded structure of $\mathcal{R}_{\ell+k}(S)$ inasmuch as that of $\mathcal{R}_{\ell+k}(S)_0~$. We finally note that, in the transitive case envisaged, the topological nature of $\mathcal{R}_{\ell+k}(S)_0$ is entirely determined by the topological nature of its isotropy group at a point. Since this group is closed in $(\Pi^0_{\ell+k}P)_y~$, we derive that $\mathcal{R}_{\ell+k}(S)_0$ is closed in $\Pi_{\ell+k}P~$. Observe, however, that the isotropy of $\mathcal{R}_{\ell+k}(S)_0$ is not necessarily the connected component, of the unit, in the isotropy of $\mathcal
{R}_{\ell+k}(S)$.

\vspace{4 mm}

\newtheorem{inf-quotient}[CorollaryCounter]{Corollary}
\begin{inf-quotient}
Let S be a structure of species E and of order $\ell$ such that $p_k\Psi(S)$ is locally of constant rank. Under these conditions,

\vspace{2 mm}

a) $\mathcal{R}_{\ell+k+h}(S)_0$ is the standard $h-th$ prolongation of $\mathcal{R}_{\ell+k}(S)_0$ i.e.,

\begin{equation*}
\mathcal{R}_{\ell+k+h}(S)_0=\Pi_{\ell+k+h}P\cap J_h\mathcal{R}_{\ell+k}(S)_0~.
\end{equation*}

\vspace{2 mm}

b) If, moreover, S verifies the equivalent conditions of the Theorem, then $\mathcal{R}_{\ell+k+h}(S)$ is the standard $h-th$ prolongation of
$\mathcal{R}_{\ell+k}(S)$. 
\end{inf-quotient}

\vspace{4 mm}

\noindent
The proof of this corollary relies on \cite{Kumpera1972} and on the following Lemma:

\vspace{2 mm}

\newtheorem{not-quotient}[LemmaCounter]{Lemma}
\begin{not-quotient}
Let $E'\longrightarrow E\longrightarrow(E'',\sigma)$ be an exact sequence of fibrations over the base space P (exact in the set theoretical sense and also in the vertical tangential sense). Then, for any integer k, the sequence of prolonged fibrations $J_kE'\longrightarrow J_kE
\longrightarrow(J_kE'',j_k\sigma)$ is also exact.
\end{not-quotient}

\vspace{4 mm}

\noindent
It is then achieved by an inductive argument on the integer \textit{k} using, at each stage, the affine structure of the kernels as well as the exactness of the sequence of linear symbols that is a consequence of the tangential exactness of the initially given sequence.

\vspace{4 mm}

\noindent
In order to prove the part (\textit{b}) of the corollary, we simply use the exactness of the sequence

\begin{equation*}
\mathcal{R}_{\ell+k}(S)~\longrightarrow~\Pi_{\ell+k}P~\xrightarrow{p_k\Phi(S)}~(J_kE,j_kS)
\end{equation*}

\vspace{2 mm}

\noindent
whose tangential exactness follows from the local constancy of the rank of $p_k\Phi(S)$ and thereafter observe that the diagram below is commutative and exact:

\begin{equation*}
\textbf{1}~\longrightarrow~J_h\mathcal{R}_{\ell+k}(S)~\longrightarrow~J_h\Pi_{\ell+k}P~\xrightarrow{J_hp_k\Phi(S)}~(J_hJ_kE,j_hj_kS)
\end{equation*}
\begin{equation*}
\hspace{50 mm}\nearrow
\end{equation*}
\begin{equation*}
\hspace{7 mm}\iota\uparrow\hspace{24 mm}\iota\uparrow\hspace{7 mm}p_hp_k\Phi(S)\hspace{12 mm}\iota\uparrow
\end{equation*}
\begin{equation*}
\hspace{23 mm}\diagup
\end{equation*}
\begin{equation*}
\textbf{1}~\longrightarrow~\mathcal{R}_{\ell+h+k}(S)~\longrightarrow~\Pi_{\ell+h+k}P~\xrightarrow{p_{h+k}\Phi(S)}~(J_{h+k}E,j_{h+k}S)
\end{equation*}
\begin{equation*}
\end{equation*}
\begin{equation*}
\uparrow\hspace{28 mm}\uparrow\hspace{32 mm}
\end{equation*}
\begin{equation*}
\end{equation*}
\begin{equation*}
\textbf{1}\hspace{28 mm}\textbf{1}\hspace{32 mm}
\end{equation*}

\vspace{4 mm}

\noindent
As for the part (\textit{a}), we simply replace, according to the Proposition 2, the previous exact sequence by

\begin{equation*}
\mathcal{R}_{\ell+k}(S)_0~\longrightarrow~\mathcal{U}~\xrightarrow{p_k\Phi(S)}~(J_kE,j_kS)
\end{equation*}

\vspace{2 mm}

\noindent
We now observe that it is essential to use hypotheses guaranteeing the appropriate differentiable structures for $\mathcal{R}_{\ell+k}(S)$ and $\mathcal{R}_{\ell+k}(S)_0$ in the lack of which the above Corollary becomes inexact, not subsisting but the inclusion

\begin{equation*}
\Pi_{\ell+h+k}P\cap J_h\mathcal{R}_{\ell+k}(S)\subset\mathcal{R}_{\ell+h+k}(S).
\end{equation*}

\vspace{2 mm}

\noindent
The above proof being rather esotheric\footnote{qualification donnée, dans les écoles des anciens philosophes, à leure doctrine secrète, réservée aux seuls initiés}, we shall transcribe it in local coordinates for the usage of the non-initiated. However, this naive transcription will be useful later. We recall that both Lie and Cartan frequently indulged into incredible calculations since they believed, presumably, that this was the first step in understanding Heaven right here from earth. Nowadays, calculations are for many a boring activity though, for others, become indispensable. Just imagine a cosmologist trying to figure out whether Einstein's constant \textit{c} is the same here in our vicinity as in, for example, \textit{Andromeda}, 60 million light years away, or in the \textit{Whirlpool Galaxy}? The same doubts arise inasmuch for $\pi$ and on how do Pythagorean circles behave in the \textit{Magellania Cloud} or on how does \textit{e} behave, does it also change + for $\times$, in \textit{Antennae}? Even at a much closer range, just beyond the neutral zone, we might ask whether the Klingon's uncertainty (undecidability) $\hbar$ is as high as for the humans or whether they are more self-confident? (\cite{Mayot1945})

\vspace{4 mm}

\begin{figure}[ht!]
\centering
\includegraphics[scale=1.7]{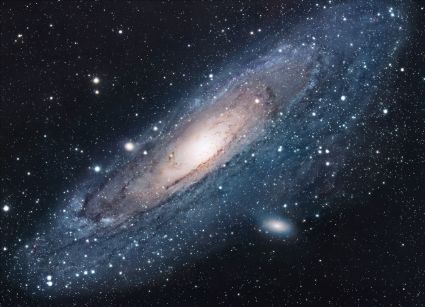}
\caption{The Universe}
\label{fig:univerise}
\end{figure}

\vspace{4 mm}

\noindent
In order to do so, let us return to the notations in the proof of the last Proposition and take a point $X\in\mathcal{R}_{\ell+k}(S)$ as well as a differentiable $\alpha-$section $\mu:\overline{\mathcal{U}}~\longrightarrow~\Pi_{\ell+k}P$ taking its values in $\mathcal{R}_{\ell+k}(S)$ and such that $\mu(y)=X~$. By the local constancy of the rank of $p_k\Phi(S)$, there exists an open neighborhood $\mathcal{U}$ of \textit{X} in $\Pi_{\ell+k+h}P$ such that $\mathcal{W}=p_k\Phi(S)(\mathcal{U})$ is a regularly embedded sub-manifold of $J_kE$ and that $p_k\Phi(S):\mathcal{U}~\longrightarrow~\mathcal{W}$ is a surmersion. We can assume, as previously, that $\overline{\mathcal{U}}=\alpha(\mathcal{U})=\alpha(\mathcal{U}\cap im~\mu)$ which implies that $\mathcal{W}\cap im~j_kS=j_kS(\overline{\mathcal{U}})$ is a regularly embedded sub-manifold of $\mathcal{W}$ . Shrinking, if necessary, the open set $\mathcal{U}$ , let us take a finite family $(f_i)$ of independent functions whose zeros define the sub-manifold $j_kS(\overline{\mathcal{U}})$ in $\mathcal{W}$ .
The local constancy of the rank of $p_k\Phi(S)$ combined with the key property $\mathcal{W}\cap im~j_kS=j_kS(\overline{\mathcal{U}})$ forces the composite functions $F_i=f_i\circ p_k\Phi(S)$ to be also independent, their zeros defining $\mathcal{R}_{\ell+k}(S)\cap \mathcal{U}~$. We shall now complete the functions $(F_i)$ to a local coordinates system by adding some functions $(G_j)$. The condition $Y=j_{\ell+h+k}\varphi(y)\in\mathcal{R}_{\ell+h+k}(S)$ means that $p_{h+k}\Phi(S)(Y)=j_{h+k}S(y)$ or, equivalently, since $p_{h+k}=(p_k)_h~$, that

\begin{equation*}
j_h(p_k\Phi(S)(j_{\ell+k}\varphi))(y)=j_h(j_kS)(y).
\end{equation*}

\vspace{4 mm}

\noindent
Translated into coordinates, this condition reads $j_h(F_i\circ j_{\ell+k}\varphi)(y)=0~$. However, still in coordinates, if we replace $j_{\ell+k}\varphi$ by the local section $\sigma$ of $\alpha:\mathcal{R}_{\ell+k}(S)~\longrightarrow~P$ whose components along the coordinates $F_i$ are null and, along the $G_j$, are equal to those of $j_{\ell+k}\varphi~$, we shall obtain the equality $j_h(j_{\ell+k}\varphi)(y)=j_h\sigma(y)$ and consequently $Y\in\Pi_{\ell+h+k}P\cap J_h\mathcal{R}_{\ell+k}(S)$. The above argument is in fact a transversality argument of $p_k\Phi(S)$ with $im~j_kS$ in a slightly more general context since transversality of the two, in the usual context, does not hold. What holds in fact is the following: There exists, in an open neighborhood $\mathcal{W}$ of each point \textit{y} belonging to $im~j_kS~$, a local foliation for which one of its leaves is an open subset of $im~j_kS$ (for example, the foliation $f_i=c_i$ of the former open set $\mathcal{W}$) and such that, if we denote by $\mathcal{W}'$ the quotient of $\mathcal{W}$ modulo the leaves and by $\rho:\mathcal{W}~\longrightarrow~\mathcal{W}'$ the projection, the composed map $\rho\circ p_k\Phi(S)$ will be of constant rank on a sufficiently small open neighborhood $\mathcal{V}$ of \textit{y}. This argument proves the claim in part (\textit{b}) of the Corollary. As for the part (\textit{a}) it will suffice to repeat the argument placing us above the open set $\mathcal{U}$.

\vspace{4 mm}

\noindent
Let us now reassume the general case where no regularity hypothesis, on $p_k\Psi(S)$, is assumed. We already remarked that the isotropies $(\textbf{R}^0_{\ell+k}(S))_y$ and $(\mathcal{R}^0_{\ell+k}S)_y$ only depend upon the jet $j_kS(y)$ and that, in particular, the isotropy of order $\ell$ only depends on the point $S(y)$. It then follows that the \textit{symbol} $(\mathfrak{g}_{\ell}S)_y$ of $\textbf{R}_{\ell}(S)$ at the point \textit{y} also depends only on $S(y)$. Observing that $\textbf{R}_{\ell+k}(S)=ker~p_k\Psi(S)$, a simple calculation (\cite{Kumpera1975}, Lemma 23.2, \cite{Kuranishi1967}, Proposition 4.3, \cite{Ruiz1977}, Proposition 9.3) will show that the symbol $(\mathfrak{g}_{\ell+k}S)_y$ of $\textbf{R}_{\ell+k}(S)$ is the $k-th$ algebraic prolongation (\textit{espace déduit}) of the symbol $(\mathfrak{g}_{\ell}S)_y$ and consequently, 

\vspace{4 mm}

\newtheorem{symbolic}[PropositionCounter]{Proposition}
\begin{symbolic}
The symbol of $\textbf{R}_{\ell+k}(S)$ at the point y only depends upon $S(y)\in E$ and this result is independent of any regularity condition requirement on the morphism $p_k\Psi(S)$. Moreover, the symbol $(\mathfrak{g}_{\ell+k}S)_y$ is the $k-th$ algebraic prolongation of $(\mathfrak{g}_{\ell}S)_y~$.
\end{symbolic}

\vspace{4 mm}

\noindent
We now examine the non-linear situation. Let $X\in\mathcal{R}_{\ell+k-1}(S)$, consider the canonical projection $\rho:\Pi_{\ell+k}P~\longrightarrow~\Pi_{\ell+k-1}P$ and define the non-linear symbol

\begin{equation*}
(g_{\ell+k}S)_X=\{Y\in\mathcal{R}_{\ell+k}(S)~|~\rho(Y)=X\}
\end{equation*}

\vspace{2 mm}

\noindent
of $\mathcal{R}_{\ell+k}(S)$ above \textit{X}. Let us show and this without any regularity hypotheses on $p_k\Phi(S)$ that $g_{\ell+k}S)_X$ is an affine sub-space of the total symbol above \textit{X} namely, the space $\{Y\in\Pi_{\ell+k}P~|~\rho(Y)=X\}$. We first observe that, for each pair $(y,z)\in P\times P$, the set

\begin{equation*}
\mathcal{R}_{\ell+k}(S)_{(y,z)}=\{X\in\mathcal{R}_{\ell+k}(S)~|~\alpha(X)=y,~\beta(X)=z\}
\end{equation*}

\vspace{2 mm}

\noindent
is a simply transitive homogeneous space of the group $\mathcal{R}^0_{\ell+k}(S)_y$ by the right action and also of the group $\mathcal{R}^0_{\ell+k}(S)_z$ by the left action. It then follows that $\mathcal{R}_{\ell+k}(S)_{(y,z)}$ is a closed and regularly embedded sub-manifold of $\Pi_{\ell+k}P_{(y,z)}$ canonically isomorphic to the left or to the right isotropy. For any $X\in\mathcal{R}_{\ell+k-1}(S)$, with $y=\alpha(X)$ and $z=\beta(X)$, the same argument shows that the symbol $(g_{\ell+k}S)_X$ is a closed and regularly embedded sub-manifold of the total symbol above \textit{X} since it is a simply transitive homogeneous space of the closed Lie group

\begin{equation*}
(g_{\ell+k}S)_y=ker~[\mathcal{R}^0_{\ell+k}(S)_y~\longrightarrow~\mathcal{R}^0_{\ell+k-1}(S)_y]
\end{equation*}

\vspace{2 mm}

\noindent
by the right action and also of the corresponding group $(g_{\ell+k}S)_z$ by the left action. Note that $(g_{\ell+k}S)_y$ is simply the symbol of $\mathcal{R}_{\ell+k}(S)$ above the unit of $\mathcal{R}_{\ell+k-1}(S)$ that identifies with \textit{y}. Let us finally show that $(g_{\ell+k}S)_X$ is a connected sub-manifold, in fact an affine sub-space of the total symbol. For this, let \textit{Z} be the projection of \textit{X} in $\mathcal{R}_{\ell}(S)$ and observe that every $Y\in(g_{\ell+k}
S)_X$ projects onto \textit{Z} by the projection $\rho_{\ell}$. We now take (\textit{at last}) local coordinate systems $(y^i)$ in an open neighborhood $\mathcal{U}$ of \textit{y} and $(z^\lambda)$ in an open neighborhood $\mathcal{V}$ of \textit{z} and consider the corresponding "jet" coordinate system $(y^i,z^{\lambda},z^{\lambda}_{\alpha})_{|\alpha|\leq r}$ on the open set $(\alpha\times\beta)^{-1}_r(\mathcal{U}\times\mathcal{V})$ of $\Pi_rP~$. Inasmuch, we also take an adapted local coordinate system $(y^i,w^{\mu})$ in an open neighborhood $\mathcal{W}$ of the point $\Phi(S)(Z)=S(y)$ in \textit{E} $(y^i=y^i\circ\pi_0)$ and there exists of course an open neighborhood $\mathcal{U}_{\ell}$ of \textit{Z} in $(\alpha\times\beta)^{-1}_{\ell}(\mathcal{U}\times\mathcal{V})$ such that $\Phi(S)(\mathcal{U}_{\ell})\subset\mathcal{W}~$. We write $\{\Phi^i,\Phi^{\mu}\}$ the components, along the coordinates $(y^i,w^{\mu})$, of the restriction $\Phi(S):\mathcal{U}_\ell~\longrightarrow~\mathcal{W}~$. Since $\Phi(S)$ is a morphism over the Identity, we infer that $\Phi^i(y^i,z^{\lambda},z^{\lambda}_{\alpha})=y^i$. We denote by $\mathcal{U}_{\ell+k}$ the inverse image, by $\rho:~\Pi_{\ell+k}P~\longrightarrow~Pi_{\ell}P~$, of the open set $\mathcal{U}_{\ell}$ and observe that the former contains 
$(g_{\ell+k}S)_X$ and is endowed with the restrictions of the coordinates $(y^i,z^{\lambda},z^{\lambda}_{\alpha})_{|\alpha|\leq\ell+k}~$. Finally, denote by $\rho^{-1}(\mathcal{W})=\mathcal{W}_k$ the open set, of $J_kE$, inverse image of $\mathcal{W}$ and equipped with the natural coordinates $(y^i,w^{\mu},w^{\mu}_{\alpha})_{|\alpha|\leq k}$ derived from $(y^i,w^{\mu})$. Then $p_k\Phi(S)(\mathcal{U}_{\ell+k})\subset\mathcal{W}_k$ and the components of the restriction $p_k\Phi(S):\mathcal{U}_{\ell+k}~\longrightarrow~\mathcal{W}_k~$, with respect to the coordinates $(y^i,w^{\mu},w^{\mu}_{\alpha})_{|\alpha|\leq k}$, are precisely the functions $\{\Phi^i,\Phi^{\mu},\partial_{\alpha}\Phi^{\mu}\}_{|\alpha|\leq k}$ where $\partial_{\alpha}$ is the \textit{total derivative} of order $|\alpha|$ with respect to the variables $(y^i)$ (iterated total derivatives in jet spaces). In particular, the symbol $(g_{\ell+k}S)_X$ is defined, in the affine space of the total symbol over \textit{X}, by the equations:

\begin{equation*}
\partial_{\alpha}\Phi^{\mu}(Y)=w^{\mu}_{\alpha}(j_kS(x))=(\partial^{\alpha} S^{\mu}/\partial y^{\alpha})(y),\hspace{5 mm}|\alpha|=k~.
\end{equation*}

\vspace{2 mm}

\noindent
Since

\begin{equation*}
\partial_{\alpha}\Phi^{\mu}(Y)=\sum_{|\beta|=\ell}~(\partial\Phi^{\mu}/\partial z^{\lambda}_{\beta})(Z)\cdot z^{\lambda}_{\alpha+\beta}~+~F_{\alpha}(X)~,
\end{equation*}

\vspace{2 mm}

\noindent
where, on the right hand side, we only detail the highest order terms, the symbol $(g_{\ell+k}S)_X$ being thereafter determined by the following linear equations with constant coefficients in the variables $z^{\lambda}_{\gamma},~|\gamma|=\ell+k~$,

\begin{equation*}
\sum_{|\beta|=\ell}~(\partial\Phi^{\mu}/\partial z^{\lambda}_{\beta})(Z)\cdot z^{\lambda}_{\alpha+\beta}=(\partial^{\alpha}S^{\mu}/\partial y^{\alpha})(y)~-~F_{\alpha}(X),
\end{equation*}

\vspace{4 mm}

\noindent
defining, as pretended, a linear affine sub-space in the space of the total symbol.

\vspace{4 mm}

\noindent
Let us now glimpse at the intrinsical aspects. With the help of the canonical identification, we can see that the abelian Lie algebra  $(\mathfrak{g}_{\ell+k}S)_z\subset\textbf{R}^0_{\ell+k}(S)_z$ is not only the Lie algebra of the group 
$(g_{\ell+k}S)_z$ but has much more impact. Recalling the results of \cite{Kumpera1975}, $\S$ 19, each element $v\in\mathfrak{g}_{\ell+k}S)_z$ determines a vector field on the total symbol space that generates a global 1-parameter group $(\varphi_t)_t$ with the property that $\varphi_1(Y)=Y+v$ is precisely the affine operation by the vector \textit{v}. The orbits of this action are all the linear affine sub-spaces of the total symbol whose direction is given by $(\mathfrak{g}_{\ell+k}S)_z$ and finally, since the sub-spaces generated, at each point, by the above vector fields are necessarily contained in $\Delta_{\ell+k}$ , the Proposition 1 will imply that an orbit, by the above (infinitesimal) affine action of $(\mathfrak{g}_{\ell+k}S)_z~$, that contains an element of $(g_{\ell+k}S)_X$ is entirely contained in $(g_{\ell+k}S)_X~$. A dimensional 
argument will also show that $dim~(g_{\ell+k}S)_X$ is equal to the dimension of 
the orbits 
and the previous tinkering (\textit{bricolage}) with local coordinates only serves 
to prove that $(g_{\ell+k}S)_X$ is also connected hence equal to an entire orbit. Furthermore, if we consider as symbol of a non-linear equation the family of all tangent spaces to the non-linear symbol $(g_{\ell+k}S)_X~$, \textit{the tangent symbol}, we perceive that this family of tangent symbols is nothing else, by the canonical identification, than $(\mathfrak{g}_{\ell+k}S)_z~$. It then follows that the tangent symbol of $\mathcal{R}_{\ell+k}(S)$ above the point \textit{X} is the $k-th$ algebraic prolongation of the tangent symbol of $\mathcal{R}_{\ell}(S)$ at 
the point \textit{Z}, hence only depends on $S(z)$. Otherwise, this last result can equally be obtained with the help of the Lemma 23.2 in \cite{Kumpera1975} or the Proposition 4.3 in \cite{Kuranishi1967} or still the Proposition 9.3 in \cite{Ruiz1977}. Let us finally observe that the groupoid structure of $\mathcal{R}_{\ell+k}(S)$ enables us  to further explicit the affine structure of 
the non-linear symbol. In fact, the argument in coordinates shows that the group $(g_{\ell+k}S)_z$ is connected, non-compact and actually homeomorphic to a numerical space. Being abelian, it canonically identifies with its Lie algebra $(\mathfrak{g}_{\ell+k}S)_z$ and the affine action of the latter on the total symbol above \textit{X} is nothing else but the left action by the abelian group $(g_{\ell+k}S)_z$. The symbol $(g_{\ell+k}S)_X$ is just one of the orbits of this action, the restricted action becoming simply transitive (without fixed points). Similarly, the right action of $(g_{\ell+k}S)_y$ on $(g_{\ell+k}S)_X$ is an affine space structure that coincides with the previous one as soon as we identify $(g_{\ell+k}S)_y$ with $(g_{\ell+k}S)_z$ by means of a conjugation via an element of $(g_{\ell+k}S)_X$.

\vspace{2 mm}

\newtheorem{affine}[PropositionCounter]{Proposition}
\begin{affine}
Without any regularity hypotheses on $p_k\Phi(S)$, the symbol of $\mathcal{R}_{\ell+k}(S)$ above any point $X\in\mathcal{R}_{\ell+k-1}(S)$ is an affine sub-space of the total symbol and the corresponding affine structure can be identified with the left action by the symbol $(g_{\ell+k}S)_z$ of the isotropy at the point $z=\beta(X)$. The tangent symbol above X is isomorphic to $(\mathfrak{g}_{\ell+k}S)_z$ and consequently only depends upon the point $S(z)\in E~$. In particular, the non-linear symbol and the tangent symbol above a unit $\textbf{e}\in\mathcal{R}_{\ell+k-1}(S)$ i.e., the symbol of the isotropy group of order $\ell+k$ at the point $y=\alpha(\textbf{e})=\beta(\textbf{e})$ as well as its Lie algebra $(\mathfrak{g}_{\ell+k}S)_y$ only depend upon $S(y)$.
\end{affine}

\section{The general problem}
Let $(E,\pi,P,p)$ be a finite prolongation space of order $\ell$ and let us now assume that \textit{E} hence consequently \textit{P} are paracompact spaces. Let us denote by $\Gamma=\Gamma(P)$ the general pseudo-group of all the local diffeomorphisms of \textit{P} and by $\mathcal{L}=\mathcal{L}(P)$ the pseudo-algebra of all the local vector fields (infinitesimal automorphisms). By prolongation of $\Gamma$ (resp. $\mathcal{L}$) to $J_kE$, we obtain the pseudo-group $\Gamma_{\ell+k}$ resp., the pseudo-algebra (pre-sheaf of Lie algebras) $\mathcal{L}_{\ell+k}$ . These are obtained by localisation of $p_k\Gamma$ resp., $p_k\mathcal{L}$ which means that we consider the set of all local finite or infinitesimal transformations of $J_kE$ that coincide locally with the prolonged transformations (and where $p_k=p_k\circ\pi$). Although $\Gamma$ and $\mathcal{L}$ are "Lie" at any order, this might fail to be true with the prolonged objects, the regularity of these being closely related to the geometry of the prolongation space \textit{E}. Nevertheless, we can still obtain much information concerning the formal equivalence problem as well as on other matters involving structures of species \textit{E} by examining closely the trajectories (orbits) of these prolonged pseudo-groups and pseudo-algebras. In fact, the "Lemma" is as follows:

\vspace{4 mm}

\textit{Two k-jets of structures of species E are equivalent (or two germs of structures of species E are equivalent up to order k) when the two jets find themselves on the same trajectory of $\Gamma_{\ell+k}~$}.

\vspace{4 mm}

\noindent
The prolongation space $p_k:J_kE~\longrightarrow~P$ being of order $\ell+k~$, we know that any local diffeomorphism $\varphi$ of \textit{P} prolongs to a local diffeomorphism $p_k\varphi$ defined by $p_k\varphi(X)= j_k(p\varphi)\cdot X$ . We thus obtain a left or right action of the groupoid $\Pi_{\ell+k}P$ on $J_kE$ though, for the time being, we only consider the left action. If $Z\cdot X=Z'\cdot X~$, then of course $Z^{-1}Z'\in \Pi^0_{\ell+k}P_X~$, the isotropy group of $\Pi_{\ell+k}P$ at the point \textit{X}, that we shall denote by $H_{\ell+k}(X)$ or simply $H(X)$. The action being differentiable, each isotropy group is a closed Lie subgroup of $\Pi^0_{\ell+k}P_y~,~y=\alpha(X)$ and the isotropies at two distinct points of the same orbit $\Omega(X)=\Pi_{\ell+k}P\cdot X$ are conjugate subgroups. Furthermore, the quotient space $(\Pi_{\ell+k}P)_y~/~H_{\ell+k}(X)$ of the classes, to the left, of $(\Pi_{\ell+k}P)_y$ that are orbits under the right action (by the source) of $H_{\ell+k}(X)$, is a differentiable fibre bundle in homogeneous spaces via the left action (by the target) of $\Pi^0_{\ell+k}P)_z~,~z=\beta(X)$, on the fibre above \textit{z}. The isotropy of this left action at the point $Z\in H(X)$ is equal to $H(Z\cdot X)$, it is obtained by the conjugation $H(Z\cdot X)=ZH(X)Z^{-1}$ and the diagram below is commutative:

\vspace{4 mm}

\begin{equation*}
(\Pi_{\ell+k}P)_y\hspace{7 mm}\xrightarrow{~\overline{q}~}\hspace{7 mm}\Omega(X)
\end{equation*}
\begin{equation*}
\downarrow\hspace{20 mm}\nearrow\mu
\end{equation*}
\begin{equation*}
(\Pi_{\ell+k}P)_y~/~H_{\ell+k}(X)\hspace{25 mm}
\end{equation*}

\vspace{4 mm}

\noindent
the arrow $\mu$ being bijective onto $\Omega(X)$ and differentiable as a mapping into $J_kE$. Let us now transport on $\Omega(X)$ and by means of $\mu$ the differentiable structure coming from the quotient and show that $\Omega(X)$ becomes a sub-manifold of $J_kE$ though not necessarily regularly embedded. To do so, it will suffice to show that $\mu$ , as a mapping into $J_kE$ , is an immersion (maximum injective rank) and this leads us to examine the infinitesimal prolongation.

\vspace{4 mm}

\noindent
The prolonged infinitesimal pseudo-algebra $\mathcal{L}_{\ell+k}$ induces, at each point, a subspace of the tangent space to $J_kE$ and consequently a distribution (field of contact elements) $\Delta_k$ on the manifold $J_kE$ that is generated by a family of vector fields stable under the bracket. Though involutive, this distribution can admit singularities. Since the infinitesimal prolongation operator $p_k$ is of order $\ell+k$, each subspace $(\Delta_k)_X$ is entirely determined by $(J_{\ell+k}TP)_y~,~y=\alpha(X)$ and, more precisely, the following sequence is exact:

\begin{equation*}
J_{\ell+k}TP~\times_P~J_kE~\xrightarrow{\lambda_k}~\Delta_k~\longrightarrow~0
\end{equation*}

\vspace{2 mm}

\noindent
Since the sheaf of germs of local sections of $J_{\ell+k}TP$ is free and of finite rank, the image sheaf, that is closed for the bracket and generates $\Delta_k~$, is also of finite type hence (\cite{Hermann1962},\cite{Turiel1976}) every $X\in J_kE$ is contained in a maximal integral sub-manifold $\omega(X)$ and verifies $T_Y\omega(X)=(\Delta_k)_Y$ for all $Y\in\omega(X)$, though the ensemble of these integral sub-manifolds does not form necessarily a regular foliation since their dimensions can vary. The space $J_kE$ admits therefore a partition by integral leaves, with eventual singularities, of $\Delta_k~$. Moreover, the leaf $\omega(X)$ is the set of points of $J_kE$ that can be joined from \textit{X} by piece-wise differentiable integral curves of $\Delta_k~$. Since $\Delta_k$ is generated by $\mathcal{L}_{\ell+k}~$, the leaf $\omega(X)$ is also the trajectory of $\mathcal{L}_{\ell+k}$ passing by \textit{X} hence, the set of all points of $J_kE$ that we can join to \textit{X} (or, for that matter, from \textit{X}) by piece-wise differentiable curves where each differentiable arc is the trajectory of a vector field belonging to $\mathcal{L}_{\ell+k}~$. We infer that $\omega(X)\subset\Omega(X)$ and, more generally, that $\Omega(X)$ is a union of integral leaves of $\Delta_k~$. Furthermore, since $\underline{J_{\ell+k}TP}$ is the Lie algebroid of the Lie groupoid $\Pi_{\ell+k}P$ and since

\vspace{4 mm}

a)the prolongation, by the target, of $\mathcal{L}$ to $\Pi_{\ell+k}P$ is infinitesimally transitive on each $\alpha-$fibre, the trajectories of the prolonged algebroid $\mathcal{L}_{\ell+k}$ being the connected components of the $\alpha-$fibres - maximal integral sub-manifolds of the trivial distribution $V\Pi_{\ell+k}P$ - as well as

\vspace{2 mm}

b) the infinitesimal action $p_k$ being derived from the finite action $p_k$ ,

\vspace{4 mm}

\noindent
we infer (always under the canonical identification) that

\vspace{2 mm}

c) the tangent map to $(\Pi_{\ell+k}P)_y~\longrightarrow~\Omega(X)\subset J_kE$ at the point \textit{Z} is equal to $\lambda_k:(J_{\ell+k}TP)_{\beta(Z)}~\longrightarrow~T_{Z\cdot X}J_kE$ and consequently its rank is equal to $dim(\Delta_k)_{Z\cdot X}~$,

\vspace{2 mm}

d) the Lie algebra $h(X)$ of $H(X)$ is equal to the kernel of the map

\begin{equation*}
\lambda_k:(J_{\ell+k}TP)_{\alpha(X)}~\longrightarrow~T_XJ_kE~,
\end{equation*}

\vspace{2 mm}

\noindent
and

\vspace{2 mm}

e) since $ZH(X)=H(Z\cdot X)Z~$, the kernel of the tangent map to $(\Pi_{\ell+k}P)_y~\longrightarrow~(\Pi_{\ell+k}P)_y~/~H_{\ell+k}(X)$ at the point \textit{Z} is equal to $h(Z\cdot X)=ker~(\lambda_k)_{\alpha(Z\cdot X)=\beta(Z)}$,

\vspace{2 mm}

\noindent
we conclude that

\begin{equation*}
(\mu_*)_{ZH(X)}:(J_{\ell+k}TP)_{\beta(Z)}~/~h(Z\cdot X)~\longrightarrow~(\Delta_k)_{Z\cdot X}
\end{equation*}

\vspace{2 mm}

\noindent
is an isomorphism hence $\Omega(X)$ is a sub-manifold of $J_kE$ for which the integral leaves of $\Delta_k$ are open sets. Since these leaves are the trajectories of $\mathcal{L}_{\ell+k}~$, they are disjoint and constitute the connected components of $\Omega(X)$. Finally, since the connected components of $\Pi_{\ell+k}P_y$ are the trajectories of the standard prolongation of $\mathcal{L}$ by the target, we see that the image of each connected component of $\Pi_{\ell+k}P_y$ by the map $\overline{q}$ is an integral leaf of $\Delta_k$ contained in $\Omega(X)$ and therefore the inverse image of a leaf is a union of connected components. In particular, the image of the connected component of the unit at the point \textit{y} is equal to $\omega(X)$.

\vspace{4 mm}

\newtheorem{orbital}[TheoremCounter]{Theorem}
\begin{orbital}
Each orbit $\Omega(X)$ of $\Gamma_{\ell+k}$ is a differentiable sub-manifold of $J_kE$ canonically isomorphic to $\Pi_{\ell+k}P_y~/~H_{\ell+k}(X)$ and invariant under $\mathcal{L}_{\ell+k}~$, the infinitesimal action being transitive. The quadruple $(\Omega(X),\alpha,P,p_k)$ is a finite prolongation space of order $\ell+k$ and the groupoid $\Pi_{\ell+k}P$ as well as the sheaf (pseudo-algebra)  $\underline{J_{\ell+k}TP}$ operate on it. The restrictions of $~\Gamma_{\ell+k}$ and $\mathcal{L}_{\ell+k}$ to $\Omega(X)$ are finite and infinitesimal pseudo-groups and pseudo-algebras of arbitrary order. For all $k\geq h$, the canonical projection $\rho_{h,k}$ transforms every $k-th$ order orbit onto an $h-th$ order orbit and thus defines a prolongation spaces morphism. The distribution $\Delta_k$ on $J_kE$ induced by $\mathcal{L}_{\ell+k}$ is involutive and locally of finite type, its maximal integral manifolds are the orbits of $\mathcal{L}_{\ell+k}$ and each orbit of the finite action has for its connected components the orbits of the infinitesimal action. The quadruple $(\omega(X),\alpha,P,p_k)$ is an infinitesimal prolongation space of order $\ell+k$ whenever \textit{P} is connected.  
\end{orbital}

\vspace{4 mm}

\noindent
The standard prolongation by the target (\cite{Kumpera1975}, $\S$ 16, item (a)) determines a canonical finite prolongation structure $(\Pi_kP,\beta,P,p^b_k)$ of finite order \textit{k} for which the $\alpha-$fibres are the trajectories\footnote{$p^s_k$ - \textit{prolongement par la source}, $p^b_k$ -  \textit{prolongement par le but}.}. For each $y\in P$, the prolongation sub-space 
$((\Pi_kP)_y,\beta,P,p^b_k)$ is transitive and the equivalence relation defined by any closed Lie sub-group $H\subset (\Pi^0_kP)_y$ is compatible with the prolongation operations. Consequently, the quotient quadruple $((\Pi_kP)_y~/~H,\beta,P,p^b_k)$ is a finite prolongation space of order \textit{k}.

\vspace{4 mm}

\newtheorem{sub-orbital}[CorollaryCounter]{Corollary}
\begin{sub-orbital}
The canonical isomorphism of the preceding theorem is an isomorphism of prolongation spaces

\begin{equation*}
((\Pi_{\ell+k}P)_y~/~H(X),\beta,P,p^b_{\ell+k})~\xrightarrow{~\mu~}~(\Omega(X),\alpha,P,p_k).
\end{equation*}
\end{sub-orbital}

\vspace{4 mm}

\noindent
Since $\Delta_k$ admits in general singularities, the space of orbits by the finite or infinitesimal actions is most often rather complicated. It can reduce to a finite number or to a discrete family of orbits (for the quotient topology), it can present itself as a regular foliation (continuous family of orbits) and, most often, the two options can appear simultaneously. The discrete orbits correspond geometrically to the existence of models (in coordinates) for the germs of structures or for their $k-$jets. Quite to the contrary, the continuous families of orbits apparently turn nonexistent the presence of models these being replaced by local deformations of non-equivalent structures, since the nature in itself of a model highlights and emphasizes the notion of rigidity (not to be confounded with the deformation of structures locally equivalent to a given model). When the orbits are discrete, the formal and local equivalence problems will have to be examined by methods specific to each case and using all the available techniques as well as the invariants. This is the case especially for the "modeled structures" (for example, modeled on a Lie pseudo-group) where we shall first try to establish the formal equivalence with the model and thereafter the local equivalence leading most often to an integrability problem. When the orbits are distributed along continuous families, it seems advantageous to appeal to the differential invariants of the action of $\Gamma$ on the space \textit{E} . These however are only speculations and the sole positive statement is the following:

\vspace{2 mm}

\textit{Two infinite jets of structures of species E are formally equivalent if and only if their $k-$jets, for any k, belong to the same $k-th$ order orbit.}

\section{The restricted problem}
Often, mainly in physics and other domains, it is important to know the equivalence not only with respect to an arbitrary transformation but also one respecting certain additional properties (\textit{e.g.}, conservation laws) and this conveys us to what we call the restricted equivalence problem with respect to the transformations of a given pseudo-group or pseudo-algebra.

\vspace{4 mm}

\noindent
Let $\Gamma$ be a pseudo-group of local transformations operating on the manifold \textit{P} and $\mathcal{L}$ the corresponding infinitesimal pseudo-algebra (sometimes called infinitesimal pseudo-group) \textit{i.e.}, the sub-presheaf of $\Gamma(TP)$ (sections) obtained by localizing as well as pasting together all the local vector fields of \textit{P} of the form $\xi=\frac{d}{dt}~\varphi_t|_{t=0}$, where $(\varphi_t)_t$ is a local one parameter family of element of $\Gamma$. We shall say that $\Gamma$ is a Lie pseudo-group of order $k_0$ if, for any $k\geq k_0$, the following properties hold:

\vspace{4 mm}

a) $J_k\Gamma$ is a closed Lie sub-groupoid of $\Pi_kP$ and the projection

\begin{equation*}
\rho:~J_{k+h}\Gamma~\longrightarrow~J_k\Gamma
\end{equation*}

\vspace{2 mm}

is a surmersion.

\vspace{2 mm}

b) $J_{k+1}\Gamma$, considered as a differential equation on the fibration

\begin{equation*}
\alpha:~\Pi_{k+1}P~\longrightarrow~P~,
\end{equation*}

\vspace{2 mm}

\noindent
is the standard prolongation of $J_k\Gamma$.

\vspace{2 mm}

c) $J_k\Gamma$ is infinitesimally complete \textit{i.e.}, the associated linear Lie equation $\textbf{R}_k=VJ_k\Gamma|P$ (\textit{P} being identified with the units of $J_k\Gamma$ and $V=\alpha-vertical$) is equal to $J_k\mathcal{L}$ .

\vspace{2 mm}

d) $\Gamma$ is complete of order $k_0$ which means that $\varphi\in\Gamma$ if and only if $j_{k_0}\varphi$ is a section (solution) of $J_{k_0}\Gamma$.

\vspace{4 mm}

\noindent
$\textbf{Remark:}$ When (a) is verified, we can easily prove that $J_{k+1}\Gamma$ is contained in the prolongation (as a differential equation) of $J_k\Gamma$ and consequently the property (b) will follow eventually at a higher order $k_0+h$ namely, when the $\delta-$cohomology of the symbols of the linear equations $\textbf{R}_k$ become $1-$acyclic. This results locally, in an open neighborhood $\mathcal{U}$ of a point in $J_{m+1}\Gamma, m=k_0+h~$, in virtue of the prolongation theorem of \textit{Cartan-Kuranishi} (\cite{Kuranishi1967}, Theorem 10.1). The invariance of the prolongation $\textit{p}(J_m\Gamma)$ by the left action of $J_{m+1}\Gamma$ shows that the open set $\mathcal{U}$ can be chosen saturated with respect to the orbits of this action. Finally, an argument based on the constancy of the characters of an exterior differential system, similar to that employed in the proof of the finiteness theorem below, shows that equality holds, at the level $m+\mu~$, in the open set $\rho^{-1}_m(\mathcal{U})$. In the next section, we provide a global proof.

\vspace{4 mm}

\noindent
The property (b) implies the corresponding property for the linear equations $\textbf{R}_k~$. Moreover, the property (d) together with (c) imply that $\mathcal{L}$ is complete of order $k_0$ and, consequently, $\mathcal{L}$ is a Lie pseudo-algebra (infinitesimal pseudo-group) of order $k_0$ since $J_k\mathcal{L}~(=\textbf{R}_k)$ is a locally trivial vector sub-bundle of $J_kTP$ and $J_{k+1}\mathcal{L}$ is the prolongation of $J_k\mathcal{L}$ in the sense of linear equations. Finally the property (c), that will be a consequence of (a), (b) and the \textit{Cartan-Kähler} theorem when the initial data is real analytic (it will also be a consequence in other situations, especially in the transitive elliptic case), serves to assure later that any orbit of the infinitesimal action is open in the corresponding orbit of the finite action and, more precisely, is a connected component.

\vspace{4 mm}

\noindent
We shall say that $\Gamma$ is \textit{transitive} when there exists $\varphi\in\Gamma$, with $\varphi(x)=y~$, whatever the points $x,y\in P$ and that it is \textit{locally trivial} when the projection $\alpha\times\beta:J_k\Gamma~\longrightarrow~P\times P$ is a submersion (we do not assume transitivity (\cite{Kumpera1971}). On account of the property (a), it will suffice to have local triviality at order $k_0~$. We shall say that $\mathcal{L}$ is \textit{transitive} or that $\Gamma$ is infinitesimally transitive when the vector sub-space induced by $\mathcal{L}$ at every point \textit{y} in \textit{P} is equal to $T_yP$. Finally, the formal transitivity is the one linked to the linear and non-linear equations $J_k\mathcal{L}$ and $J_k\Gamma$ and coincides entirely with the transitivity above.

\vspace{4 mm}

\noindent
Let $(E,\pi,P,p)$ be a finite or infinitesimal prolongation space and $\Gamma$, resp. $\mathcal{L}$, a finite or infinitesimal Lie pseudo-group (pseudo-algebra) of order $k_0$ operating on \textit{P}. The Definition 1 can be re-written by replacing the general pseudo-group and pseudo-algebra of all local finite or infinitesimal transformations by the more specific data $\Gamma$ and $\mathcal{L}~$. In this context, we can re-write essentially all of the section 2 by replacing $\Pi_{\ell+k}P$ and $J_{\ell+k}TP$ by $J_{\ell+k}\Gamma$ and $J_{\ell+k}\mathcal{L}$ as soon as $\ell+k\geq k_0~$. We can also transcribe the considerations of the previous section where we shall replace $\Gamma_{\ell+k}$, resp. $\mathcal{L}_{\ell+k}$, by the prolongations of $\Gamma$, resp. $\mathcal{L}~$, the general equivalence problem by the restricted one and, in the Theorem 2, the groupoid $\Pi_{\ell+k}P$ by $J_{\ell+k}\Gamma$. Inasmuch, we can rewrite the sections 3 and 4 in the restricted context though certain parts and especially those concerning the morphisms $\Phi$ and $\Psi$ require some additional considerations. We shall return to this in a later section.

\vspace{4 mm}

\noindent
Still in a wider context, we can define finite and infinitesimal prolongation spaces \textit{relative} to given finite or infinitesimal pseudo-groups of transformations. In other terms, the prolongation operations are only defined for the elements of the pseudo-group or pseudo-algebra envisaged. A most relevant example is provided by the Cartan \textit{normal} prolongation spaces associated to given pseudo-groups and their quotient spaces. In replacing the general pseudo-group by a given one we can still argue as in the previous sections though, of course, we shall not forget the inequality $\ell+k\geq k_0~$.

\vspace{4 mm}

\noindent
The Lie pseudo-group $\Gamma$ is said to be analytic when the manifolds  $J_{\ell+k}\Gamma,~k\geq k_0~$, are analytic sub-groupoids of $\Pi_{\ell+k}P$ (supposing of course that $\pi:E~\longrightarrow~P$ is an analytic fibration). Inasmuch, the Lie pseudo-algebra $\mathcal{L}$ is said to be analytic when the linear equations $J_k\mathcal{L}$ are analytic vector sub-bundles of $J_kTP$. Clearly, the analiticity of $\Gamma$ implies that of $\mathcal{L}$ and the converse is also true since the differentiable structure of $J_k\Gamma$ is entirely determined, in a neighborhood of the units hence everywhere, by the structure of $\textbf{R}_k~$.

\section{The role of the differential invariants - finiteness theorems}

The interesting situation from the point of view of the differential invariants is that of continuous families of orbits. We therefore assume, for the time being, that there exists an integer $k_1$ such that, for $k\geq k_1$, the orbits of the action of $\Gamma_{\ell+k}~$ on $J_kE$ or rather those of the infinitesimal action of $\mathcal{L}_{\ell+k}~$ are distributed along a regular foliation \textit{i.e.}, the integrable distribution $\Delta_k$ is locally of constant dimension. A first integral of $\Delta_k$ (a function that is locally constant on each integral leaf of $\Delta_k$ or, equivalently, a function whose differential $df$ vanishes on $\Delta_k$) will be called a \textit{differential invariant of order k} of the Lie pseudo-group $\Gamma$, resp. of the pseudo-algebra $\mathcal{L}~$, and relative to the prolongation space \textit{E}. Since $\Delta_k$ is assumed to be regular there exists, in a neighborhood of each point in $J_kE~$, a fundamental system of independent differential invariants their number (rank) being equal to the
co-dimension of $\Delta_k~$. On the other hand, with the aid of the formal derivatives (total derivatives) in the jet manifolds, it is possible to ascend (lift), in a non-trivial way, differential invariants defined on $J_kE$ to new differential invariants defined on $J_{k+1}E~$. Lie's finiteness theorem for the differential invariants states essentially that the invariants of any order are generated by those of a certain finite order together with all their successive formal derivatives. The mechanism involving the differential invariants presumes of course certain regularity hypotheses as well as specific technicalities that eventually will lead us to the Fundamental Theorem of Sophus Lie (\cite{Lie1884}) and we shall try to describe these in the most succinct manner by referring as much as possible to \cite{Kumpera1975}. Since the specific case of prolongation spaces and of the formal equivalence of structures, our main concern, isn't but a special case of the general problem discussed in the previous reference, it is possible to simplify two of the hypotheses and the form under which we state the Lie Theorem for its applications in the equivalence problem. We shall in fact provide a much more precise statement than the one claimed in the Theorem 23.6 (\textit{loc.cit.}).

\vspace{4 mm}

\noindent
Let us first remark that the problem in \cite{Kumpera1975} consists in taking an arbitrary fibration (surmersion) $~\pi:P~\longrightarrow~M~$ together with a sheaf $\mathcal{L}$ (Lie sheaf) of vector fields on \textit{P} and thereafter study the differential invariants of $\mathcal{L}$ in the realm of the standard prolongation spaces $J_kP$ above \textit{P}. Here, quite to the contrary, we are given an infinitesimal Lie pseudo-algebra $\mathcal{L}$ of order $k_0$
on the manifold \textit{P}, an infinitesimal prolongation space $(E,\pi,P,p)$ of finite order $\ell$ and study the differential invariants of the infinitesimal action of $\mathcal{L}$ on the prolongation spaces $J_kE~$. We can re-conduce our considerations to the above mentioned context by simply considering the prolonged sheaf $\mathcal{L}_{\ell}=p\mathcal{L}$ defined on the space \textit{E} and study the differential invariants of $\mathcal{L}_\ell$ by the techniques and methods found in \cite{Kumpera1975}. However, the present methods are far more reaching and accurate.

\vspace{4 mm}

\noindent
We examine initially the hypothesis $H_1$ (\textit{loc.cit.}, pg.363) or, by preference, the weaker hypothesis on pg.378.

\vspace{2 mm}

$H'_{1,Y}$: There exists an integer $k_1$ such that, for any $k\geq k_1~$, the fibre space $(\tilde{L}_V)_k$ with base space $J_{k+1}E$ is of constant rank in the neighborhood of each point $Y_{k+1}~$.

\vspace{2 mm}

\noindent
Let us recall (see \cite{Kumpera1975} for the notations) that $(\tilde{L}_V)_k\subset J_{k+1}P~\times_P~\tilde{J}_kVE$ and that this fibre space is the image of $\tilde{L}_k=\tilde{J}_k\mathcal{L}$, $\mathcal{L}$ being a Lie sheaf over \textit{P}, by the verticalisation operation described in terms of the exact sequence (22.29) in \cite{Kumpera1975}. However, in the present case\footnote{The "tilde" notation refers to the composite fibration $TE~\longrightarrow~E~\longrightarrow~P$ and where $TE$ is also replaced by $VE$}, we start with an infinitesimal Lie pseudo-algebra $\mathcal{L}$ of order $k_0$ defined on \textit{P}, consider its prolongation $\mathcal{L}_\ell=p\mathcal{L}$ to \textit{E} and $\tilde{L}_k$ becomes $\tilde{J}_k(\mathcal{L}_\ell)$. Under these conditions, $(\tilde{L}_V)_k$ is the image of $J_{k+1}E~\times_P~J_{\ell+k}\mathcal{L}$ by the mapping

\begin{equation*}
J_{k+1}E~\times_P~J_{\ell+k}TP~\longrightarrow~J_{k+1}E~\times_P~\tilde{J}_kVE
\end{equation*}

\vspace{2 mm}

\noindent
defined by

\begin{equation*}
(j_{k+1}\sigma(y),j_{\ell+k}\xi(y))~\longmapsto~j_k[(p\xi)\circ\sigma-(T\sigma\circ\xi)](y).
\end{equation*}

\vspace{2 mm}

\noindent
Let us next consider the exact sequence

\begin{equation*}
0~\longrightarrow~\mathcal{N}_{\ell+k}~\longrightarrow~J_{k+1}E~\times_P~J_{\ell+k}\mathcal{L}~\longrightarrow~(\tilde{L}_V)_k~\longrightarrow~0
\end{equation*}

\vspace{2 mm}

\noindent
where $\mathcal{N}_{\ell+k}$ denotes the kernel. We thus see that the regularity of $(\tilde{L}_V)_k$ can be replaced, when $\ell+k\geq k_0$, by that of $\mathcal{N}_{\ell+k}$ for which the defining equation is given by

\begin{equation*}
j_k[(p\xi)\circ\sigma-(T\sigma\circ\xi)](x)=0.
\end{equation*}

\vspace{2 mm}

\noindent
This equation can be envisaged as a linear differential equation of order $\ell+k$ in $J_{\ell+k}\mathcal{L}$ (or as well in $J_{\ell+k}TP$) whose coefficients depend on the parameters in $J_{k+1}E~$. Consequently, the regularity of this equation in the neighborhood of a jet $Z_{k+1}\in J_{k+1}E$, is closely related to the geometry of the prolongation space \textit{E} in the neighborhood of $\beta(Z_{k+1})$.

\vspace{2 mm}

\noindent
As for the other two hypotheses on the pg.363, we can partly weaken $H_{2,X}$ by taking into account that $J_{\ell+k}\mathcal{L}$ is a Lie equation hence the distribution $\Delta_k$ automatically satisfies the involutivity condition. However, we shall be forced to strengthen the part concerning regularity. Inasmuch, we shall strengthen the point-wise hypothesis $H_{3,X}$ by a local condition, its most efficacious verification criterion being provided by the Proposition 25.4 in \cite{Kumpera1975} on account of its Corollary. We therefore consider the following hypotheses:

\vspace{4 mm}

$H_1:$ For any $k\geq k_1~(\geq k_0-\ell),$ the vector bundle $\mathcal{N}_{\ell+k}$ has constant rank in a neighborhood of $Z_{k+1}$ and we denote by $k_1(Z)$ the integer where-after $(\Delta_{k-1,k})_{Z_k}$ (the kernel) becomes involutive.

\vspace{4 mm}

$H_2:$ There exists a family $(\mathcal{U}_k)_{k\geq k_2(Z}$ of open neighborhoods of the jets $Z_k$ such that $~\rho_{k,k+h}:\mathcal{U}_{k+h}~\longrightarrow~\mathcal{U}_k~$ is a fibration and $\Delta_k$ has constant dimension on $\mathcal{U}_k$ .

\vspace{4 mm}

$H_3:$ $\beta(\textbf{R}_{k_3(Z)}(\mathcal{L}))_{Z'})=T_{\beta(Z')}P$ for all $Z'\in\mathcal{U}_{k_3}(Z)$ and for some integer $k_3\geq k_2$ .

\vspace{4 mm}

\noindent
The last hypothesis assures the existence of a local basis of \textit{admissible} formal derivations of order $k_3$ centered around \textit{Z}, admissible meaning that such derivations transform differential invariants into differential invariants. Before stating the desired theorem, let us examine a little closer the above hypotheses in order to better discern their meaning.

\vspace{5 mm}

\noindent
1. Let $(\Delta_{k-1,k})_{Z'}$ be the kernel of $~T\rho_{k-1,k}:
\Delta_k~\longrightarrow~\Delta_{k-1}~$ at the point $Z'\in J_{k+1}E$ (this mapping being always surjective). The first hypothesis serves to prove that there is an order $k'$ such that the kernel $(\Delta_{k,k+1})_{Z_{k+1}}~$, $k\geq k'$, is contained, by means of the canonical identification, in the algebraic prolongation of $(\Delta_{k-1,k})_{Z_k}$ and consequently becomes equal to it from an order $k''$ onwards or, in other terms, the Spencer $\delta-$complex constructed with the kernels $(\Delta_{k-1,k})_{Z_k}$ becomes $1-$acyclic for $k\geq k''$. Likewise, it will become involutive beginning with an integer that we shall denote by $k_1(Z)$. This hypothesis alone enables us to prove the \textit{asymptotic stability} result (\cite{Kumpera1975}, Theorems 22.1 and 23.1) that in turn and with the aid of the hypotheses $H_{2,Z}$ and $H_{3,Z}$, leads to the Lie Theorem (\textit{loc.cit.}, Theorem 23.6).

\vspace{4 mm}

\noindent
2. The hypothesis $H_2$ assures a sufficient number \textit{i.e.}, a complete set of $k-th$ order differential invariants defined on the open set $\mathcal{U}_k~$.

\vspace{4 mm}

\noindent
3. The hypothesis $H_3$ enables us to obtain a sufficient number of $(k+1)-st$ order differential invariants by taking admissible formal derivations of the $k-th$ order differential invariants (and, of course, lifting also the latter up to order $k+1$).

\vspace{4 mm}

\noindent
4. One shows that the finiteness property of the differential invariants takes place at the stage $k~\rightsquigarrow~k+1$ (\textit{i.e.}, for the germs of invariants at the points $Z_k$ and $Z_{k+1}$ respectively) if and only if $(\Delta_{k,k+1})_{Z_{k+1}}$ is the algebraic prolongation of $(\Delta_{k-1,k})_{Z_k}$ (\cite{Kumpera1975}, Lemmas 23.3 and 23.5).

\vspace{4 mm}

\noindent
We next remark that the three hypotheses, \textit{per se} independent, are not strictly necessary to prove the desired results. In fact, the hypothesis $H_{2,X}$ underlying the theorem 23.8 in \cite{Kumpera1975} is considerably weaker than $H_2$ . However, the asymptotic stability, consequence of $H'_{1,Z}~$, joint to $H_{3,Z}$ imply the local regularity and the integrability of the distribution $\Delta_k~$, for $k>k_2$ , in view of the Corollary 5 (\textit{loc.cit.}, pg.377, conditions (I) and (II)). Viewed from another angle, we note that solely the hypotheses $H_2$ and $H_3$ will, in virtue of the lemma 23.3 in \cite{Kumpera1975}, imply that $(\Delta_{k,k+1})_{Z_{k+1}}\subset p(\Delta_{k-1,k})_{Z_k}$ and we thus obtain the asymptotic stability of the kernels starting from a certain integer $k''$. These remarks simply show that the usage of the above three hypotheses admits a certain flexibility, the appropriate choices being conditioned to the results looked for.

\vspace{4 mm}

\noindent
At present we choose $H_2$ and $H_3$ as underlying hypotheses and fix the order $\mu=k_1(Z)$ where after the symbols $(\Delta_{k-1,k})_{Z_k}$ become involutive (the hypothesis $H_1$ only reappearing later when the regularity of the $\Delta_k$ becomes apparent). We can further assume that $k_2(Z)~<~k_1(Z)~$. The hypothesis $H_2$ implies that the kernels $\Delta_{\mu,\mu+1}$ and $\Delta_{\mu-1,\mu}$ are of constant dimension in $\mathcal{U}_{\mu+1}$ and $\mathcal{U}_\mu$ respectively, and further $(\Delta_{\mu-1,\mu})_{Z_\mu}$ is involutive, $(\Delta_{\mu,\mu+1})_{Z_{\mu+1}}$ being its algebraic prolongation. According to the Theorem 23.6 (\textit{loc.cit.}), there exists an open neighborhood $\mathcal{U}_{\mu+1}$ of $Z_{\mu+1}$ such that $~(\Delta_{\mu,\mu+1})_{Z'_{\mu+1}}=p(\Delta_{\mu-1,\mu})_{Z'_{\mu}}~$ for all $Z'_{\mu+1}\in\mathcal{U}_{\mu+1}$ and consequently the finiteness property of the differential invariants is verified at the step $~\mathcal{U}_{\mu}~\rightsquigarrow~\mathcal{U}_{\mu+1}~$, $\mathcal{U}_{\mu}=\rho(\mathcal{U}_{\mu+1})$. Let us next observe that the \textit{characters} $\tau_i$ of $~(\Delta_{\mu-1,\mu})_{Z'_{\mu}}~,~Z'_{\mu}\in\mathcal{U}_{\mu}~$, are lower semi-continuous. The dimensions of $\Delta_{\mu-1,\mu}$ and $\Delta_{\mu,\mu+1}$ being constant, the characterization of the involutivity (\cite{Kumpera1975}, $\S$ 24, property 8, \cite{Kuranishi1967}, proposition 6.1) implies the existence of an open neighborhood $\mathcal{W}_{\mu}$ of $Z_{\mu}$ in which the kernels $(\Delta_{\mu-1,\mu})_{Z'_{\mu}}~$, $Z'_{\mu}\in\mathcal{W}_{\mu}~$ are all involutive, the characters $\tau_i$ remaining constant. Let us denote by $\mathcal{W}_{\mu+1}$ the inverse image of $\mathcal{W}_{\mu}$ with respect to the projection $~\rho:\mathcal{U}_{\mu+1}~\longrightarrow~\mathcal{U}_{\mu}~$ and, similarly, define $\mathcal{W}_{\mu+h}$ considering $~\rho:\mathcal{U}_{\mu+h}~\longrightarrow~\mathcal{U}_{\mu}~$. Furthermore, denote by $(\Delta'_{\mu+2})_{Z'}$ the sub-space of $T_{Z'}J_{\mu+2}E~$, ${Z'}\in\mathcal{W}_{\mu+2}~$, defined by
\begin{equation*}
(\Delta'_{\mu+2})_{Z'}=ker_{Z'} \{\rho^*_{\mu+1,\mu+2}df, \partial_\varphi  df~|~Z'_{\mu}\in\mathcal{W}_{\mu},~ f\in(\mathfrak{I}_{\mu+1})_{Z''},
\end{equation*}
\begin{equation*}
\varphi\in\mathcal{R}_{\mu+1}(\mathcal{L})_{Z''},~Z''=\rho_{\mu+1,\mu+2}Z'\},
\end{equation*}
where $\mathfrak{I}_{\mu+1}$ denotes the algebra of all differential invariants of order $\mu+1~$. Since the elements of $\mathcal{R}_{\mu+1}(\mathcal{L})_{Z''}$ are admissible, the inclusion $~\Delta'_{\mu+2}\supset\Delta_{\mu+2}~$ holds and the lemma 23.3 in \cite{Kumpera1975} shows furthermore that

\begin{equation*}
dim~(\Delta'_{\mu+2})_{Z'}=dim~(\Delta_{\mu+1})_{Z''}+dim~p(\Delta_{\mu,\mu+1})_{Z''}~.
\end{equation*}

\vspace{2 mm}

\noindent
However, $(\Delta_{\mu,\mu+1})_{Z''}$, $Z''\in\mathcal{W}_{\mu+1}$, is involutive it being the prolongation of an involutive space and the property 9 in $\S$ 24 of \cite{Kumpera1975} or else, the Proposition 9.4 in \cite{Kuranishi1967} shows, in view of the constancy of the characters $\tau_i~$, that $~dim~p(\Delta_{\mu,\mu+1})_{Z''}~$ is constant in $~\mathcal{W}_{\mu+1}~$, the characters of these prolonged spaces being also constant. Returning to the point $Z_{\mu+2}~$, we perceive that this dimension is equal to $~dim~(\Delta_{\mu+1,\mu+2})_{Z_{\mu+2}}$ and consequently that $~dim~(\Delta'_{\mu+2})_{Z'}=dim~(\Delta_{\mu+2})_{Z'_{\mu+2}}~$. Furthermore, this entails, in virtue of the constancy of the dimensions of $\Delta_{\mu+2}~$, that $~(\Delta'_{\mu+2})_{Z'}=(\Delta_{\mu+2})_{Z'}$ for all $Z'\in\mathcal{W}_{\mu+2}~$. We thus infer that the finiteness property for the differential invariants is verified at the step $~\mathcal{W}_{\mu+1}~\rightsquigarrow~\mathcal{W}_{\mu+2}~$. An inductive argument will finally prove, based on the constancy of the characters, that the finiteness property is verified at the stage $~\mathcal{W}_{\mu+h}~\rightsquigarrow~\mathcal{W}_{\mu+h+1}$ since the involutivity as well as the constancy of the characters is preserved by prolongation.

\vspace{2 mm}

\noindent
Let us observe that the involutivity property of the kernels $(\Delta_{\mu-1,~\mu})_{Z'_\mu}~$ together with the regularity of the $\Delta_k~$,  $k\geq\mu-1~$, on the open sets $\mathcal{U}_k~$, that "fibrate" one upon the other, serve uniquely to ensure the existence of a family of open neighborhoods $\mathcal{W}_k$ of $Z_k~$, fibering one above the other in such a way that, along every element $Z'\in~proj~lim~\mathcal{W}_k~$, the consecutive kernels of the $\Delta_k$ constitute a $1-$acyclic Spencer $\delta-$com-plex. In the applications, this local $1-$acyclicity property, sole to assure the finiteness mechanism of the differential invariants, might be verified long before the involutivity. The argument as well as the aims of the above discussion are somewhat quite the opposite of what has been looked for in the $\S$ 24 of \cite{Kumpera1975} where the problem posed was the regularity of the trajectories.

\vspace{4 mm}

\newtheorem{finiteness}[TheoremCounter]{Theorem (of finiteness)}
\begin{finiteness}
Let $(E,\pi,P,p)$ be an infinitesimal prolongation space, $Z\in J_{\infty}E$ an infinite jet of a structure of species E and $\mathcal{L}$ an infinitesimal pseudo-algebra (Lie pseudo-algebra) operating on P. Assuming that the hypotheses $H_2$ and $H_3$ are satisfied at the point Z, we write $\mu=k_1$ and take a family of n $(=dim~P)$ local sections $\varphi_i$ of $\mathcal{R}_{\mu+1}(\mathcal{L})$, defined in a neighborhood of $Z_{\mu+1}~$, such that $\{\beta\varphi_i(Z_{\mu+1})\}$ generates the tangent space $T_yP~,~y=\alpha(Z)$ \textit{i.e.}, the family $\{\varphi_i\}$ is a local basis of admissible formal derivations in the neighborhood of $Z_{\mu+1}~$. Under these conditions, there exists a family $(\mathcal{W}_k)_{k\geq\mu}~$, each $\mathcal{W}_k$ being an open neighborhood of $Z_k~$, such that:

\vspace{4 mm}

\hspace{5 mm}$i)\hspace{2 mm}\rho_{k,k+h}:~\mathcal{W}_{k+h}~\longrightarrow~\mathcal{W}_k$ is a fibration.

\vspace{3 mm}

\hspace{4 mm}$ii)\hspace{2 mm}(\Delta_{k,k+1})_{Z'_{k+1}}=p(\Delta_{k-1,k})_{Z'_k}~,\hspace{2 mm}Z'_{k+1}\in\mathcal{W}_{k+1}$.

\vspace{3 mm}

\hspace{3 mm}$iii)\hspace{2 mm}(\Delta_k)_{Z'_k}=ker_{Z'_k}\{df~|~f\in\mathfrak{I}\},\hspace{2 mm}Z'_k\in\mathcal{W}_k$.

\vspace{3 mm}

\hspace{4 mm}$iv)\hspace{2 mm}(\Delta_{k+1})_{Z'_{k+1}}\!=ker_{Z'_{k+1}}\{\rho^*_{k,k+1}df,~\partial{\varphi_i}df~|~f\in\mathfrak{I}_k,~1\leq i\leq n\}$,

\vspace{3 mm}

\hspace{11 mm}$Z'_{k+1}\in\mathcal{W}_{k+1}$.

\vspace{3 mm}

\hspace{5 mm}$v)\hspace{2 mm}(\Delta_{k+1})_{Z'_{k+1}}\!=ker_{Z'_{k+1}}\{\rho^*_{k,k+1}df,~\partial{\varphi_i}df~|~f\in\mathfrak{I}_k$,

\vspace{3 mm}

\hspace{11 mm}$\varphi_i\in\mathcal{R}_{\mu+1}(\mathcal{L})_{Z'_{\mu+1}}\},\hspace{2 mm}Z'_{k+1}\in\mathcal{W}_{k+1}$.

\vspace{3 mm}

\hspace{4 mm}$vi)\hspace{2 mm}\{\rho^*_{k,k+1}df,~\partial{\varphi_i}df~|~f\in\mathfrak{I}_k~,~1\leq i\leq n\}_{Z'_{k+1}}$

\vspace{3 mm}

\hspace{11 mm}\textit{generates} $(d\mathfrak{I}_{k+1})_{Z'_{k+1}}),\hspace{2 mm}Z'_{k+1}\in\mathcal{W}_{k+1}$.
\end{finiteness}

\vspace{4 mm}

\noindent
The interest of the finiteness theorem for the equivalence of structures is due to the fact that it enables us to translate the equivalence by \textit{only} a finite number of conditions (equality of the values taken by a finite number of differential invariants). We terminate here this awfully long "preamble" and will retake the effective study of the equivalence problem in part II of this paper where diverse \textit{mises en scène} shall be examined. As a last word, we should say that all the previous discussion can be carried out in the context of prolongation spaces \textit{relative} to given finite Lie pseudo-groups or infinitesimal Lie pseudo-algebras and also it is worthwhile to recall that Sophus Lie provided some of the most remarkable contributions. Surprisingly, the formula 25.5, concerning the bracket of formal and holonomic derivations (\cite{Kumpera1975}) is already written in his work \cite{Lie1884} (see also \cite{Kumpera1967}, \cite{Molino1972}).

\bibliographystyle{plain}
\bibliography{references}

\end{document}